\documentclass[11pt]{article}
\usepackage[a4paper,margin=1in]{geometry}
\usepackage{setspace}

\RequirePackage{amsthm,amsmath,amsfonts,amssymb}
\RequirePackage{mathtools,bm,aliascnt}
\RequirePackage{booktabs,tabularx,longtable,array}
\RequirePackage{algorithm}
\RequirePackage{algpseudocode}
\RequirePackage{tikz}
\RequirePackage{multirow,float,comment,placeins}
\RequirePackage{xcolor}
\RequirePackage[numbers,sort&compress]{natbib}
\RequirePackage[colorlinks=true,allcolors=blue]{hyperref}

\definecolor{revieworange}{RGB}{204,102,0}

\newtheorem{theorem}{Theorem}[section]
\newaliascnt{proposition}{theorem}
\newtheorem{proposition}[proposition]{Proposition}
\aliascntresetthe{proposition}
\newaliascnt{lemma}{theorem}

\aliascntresetthe{lemma}
\newaliascnt{corollary}{theorem}
\newtheorem{corollary}[corollary]{Corollary}
\aliascntresetthe{corollary}
\theoremstyle{definition}
\newaliascnt{assumption}{theorem}
\newtheorem{assumption}[assumption]{Assumption}
\aliascntresetthe{assumption}
\newaliascnt{definition}{theorem}

\aliascntresetthe{definition}
\newaliascnt{example}{theorem}

\aliascntresetthe{example}
\newaliascnt{remark}{theorem}
\newtheorem{remark}[remark]{Remark}
\aliascntresetthe{remark}

\newcommand{\R}{\mathbb R}
\newcommand{\N}{\mathbb N}
\newcommand{\E}{\mathbb E}
\newcommand{\Pp}{\mathbb P}
\newcommand{\cH}{\mathcal H}
\newcommand{\cX}{\mathcal X}
\newcommand{\cR}{\mathcal R}
\newcommand{\cA}{\mathcal A}

\newcommand{\Cov}{\operatorname{Cov}}

\newcommand{\tr}{\operatorname{tr}}
\newcommand{\rank}{\operatorname{rank}}

\newcommand{\dd}{\mathrm d}

\newcommand{\eps}{\varepsilon}
\newcommand{\starone}{\star_1^1}

\newcommand{\dK}{d_{\mathrm K}}
\newcommand{\dL}{d_{\mathrm L}}

\title{Approximation Theorems for High-Dimensional Canonical \textit{U}-Statistics:
Gaussian Chaos and Phase Transition}
\author{Leheng Cai\thanks{Department of Statistics and Data Science, Tsinghua University.}
\and Qirui Hu\thanks{\begin{tabular}[t]{@{}l@{}}
School of Statistics and Data Science, Shanghai University of Finance and Economics.\\
Department of Mathematics, Ruhr-Universit\"at Bochum\\
Corresponding author. Email: \href{mailto:huqirui@mail.shufe.edu.cn}{huqirui@mail.shufe.edu.cn}
\end{tabular}}}
\date{}
\begin{document}
\maketitle
\begin{abstract}
We study simultaneous inference for maxima of  canonical
order-two $U$-statistics in high dimension.  Degeneracy makes quadratic fluctuations leading, so ordinary Gaussian calibration
can fail even after exact variance normalization.  We show that the
appropriate general target is a joint signed Gaussian quadratic chaos and
establish a  general approximation result that permits indefinite kernels. The general anti-concentration bound is too crude for  high-dimensional inference, and we obtain  sharper bounds under additional spectral structure.
We also identify a  phase transition from a non-Gaussian signed-chaos maximum to its covariance-matched Gaussian counterpart driven by the effective rank.
 For feasible inference, we propose a Gaussian multiplier bootstrap that avoids estimating eigensystems, and establish its validity.  Two applications and extensive numerical simulations further illustrate the scope and practical performance of the proposed framework.
\end{abstract}



\section{Introduction}\label{sec:introduction}

Simultaneous screening of quadratic signals leads naturally to maxima of
canonical $U$-statistics.  In a high-dimensional covariance graph, for
example, the unbiased estimator of a squared edge signal is an order-two
$U$-statistic after centering the marginal covariance score.  Simultaneous
kernel two-sample testing has the same structure because the kernel of the unbiased estimator of squared MMD is canonical under the two-sample null.  In both settings,
standardization alone does not justify a normal approximation; see Section \ref{sec:applications} for details.

The obstacle is the vanishing first Hoeffding projection: the leading term
is quadratic rather than linear.  For a fixed square-integrable canonical kernel, the appropriately normalized \(U\)-statistic converges to a weighted sum of centered chi-square variables
\cite{Hoeffding1948,Neuhaus1977,ArconesGine1993}.  In the rank-one case the
variance-standardized statistic converges to $(Z^2-1)/\sqrt2$, where $Z$ is a standard normal random variable,  so exact
standardization does not restore normality.  High-dimensional Gaussian
approximations for sums of independent vectors therefore do not directly
solve the problem, despite their logarithmic dependence on the ambient
dimension \cite{CCK2013,CCK2015,CCK2017,CCKK2022,CCK2023}.

Existing high-dimensional $U$-statistic theory mainly treats nondegenerate
complete statistics, randomized incomplete statistics, or regimes in which
the quadratic statistic itself Gaussianizes
\cite{Chen2018,ChenKato2019,HuangEtAl2023,ImaiKoike2025}.  De Jong and
fourth-moment theories identify ordinary Gaussian endpoints under
vanishing-influence or contraction conditions
\cite{deJong1987,NourdinPeccatiReinert2010,DoblerPeccati2017,
DoblerPeccati2019,Koike2023}.  Those results are indispensable in diffuse
spectral regimes, but they do not supply the correct calibration when some
coordinates retain a few dominant  eigenvalues.

In this paper, we  retain the signed Karhunen--Lo\'eve expansion of each kernel separately, without requiring a common eigenbasis or simultaneous diagonalization. A single isonormal Gaussian process on the closed span of all eigenfunctions generates the joint chaos target and preserves the cross-kernel covariance structure. This construction accommodates signed spectra, finite or infinite rank, and singular dependence  among coordinates.

To establish the  theoretical results, we avoid coordinatewise spectral truncation. We   approximate all normalized kernels simultaneously by conditional expectations on a common refining sequence of finite \(\sigma\)-fields.
The approximants preserve symmetry and canonicality and share a common finite-dimensional score representation, without requiring positivity.
This permits a blockwise Lindeberg replacement, with scalar and Hilbert-space Rosenthal inequalities controlling the replacement error uniformly over the approximation level. An exact Gaussian pair-sum identity then identifies the Gaussian replacement with the corresponding quadratic-chaos target. Uniformity over the partition level allows us to pass to arbitrary rank without imposing an eigenvalue-tail condition.


Our general result first yields a L\'evy approximation, which provides a weak distributional comparison without requiring anti-concentration. Upgrading this result to a Kolmogorov approximation requires control of the probability that the chaos maximum falls in a short interval. Without additional spectral structure, a coordinatewise small-ball argument gives an anti-concentration bound with linear dependence on \(p\), and a negative rank-one example shows that this dependence is sharp in general. Sharper Kolmogorov conclusions become available under additional spectral structure, either for positive spectra with finite rank or regular spectral decay, or for signed spectra with sufficiently large effective rank.

The Gaussian chaos approximation target exhibits a spectral phase transition. The fourth-order effective rank, defined in Section~\ref{sec:setup}, quantifies the diffuseness of the  spectrum. Using a Malliavin--Stein comparison for maxima of Wiener functionals \cite{Koike2019Wiener}, we show that the signed-chaos maximum approaches its covariance-matched Gaussian counterpart once the minimum effective rank dominates the relevant logarithmic factor.


For practical implementation, we develop a  Gaussian multiplier bootstrap for positive-semidefinite kernels that avoids estimating the oracle eigensystem. The procedure is constructed from centered kernel Gram matrices and a common Gaussian multiplier vector, and does not require knowing in advance whether the relevant approximation is a non-Gaussian quadratic chaos or an ordinary Gaussian law. The same calibration therefore applies across spectral regimes: in low-effective-rank settings it reproduces the quadratic-chaos behavior, while in sufficiently diffuse regimes the corresponding chaos  Gaussianizes. Its   validity is established through a comparison between empirical and population block covariance operators, which reduces the bootstrap error to a single uniform covariance discrepancy that can be verified under different settings.

Our two applications place the theory in established high-dimensional testing problems. In covariance-structure testing, maxima of sample correlations and related coherence statistics are classical tools \cite{CaiJiang2011,CaiLiuXia2013}. After squaring and centering, the corresponding statistics are rank-one canonical \(U\)-statistics, whose natural coordinatewise limit is a centered chi-square law even under Gaussian sampling. In kernel two-sample testing, the MMD statistics estimate  distributional discrepancies \cite{GrettonEtAl2012}. While \citet{GaoShao2023} study Gaussian approximation for a single studentized MMD in high data dimension, our challenge arises from many coordinate- or blockwise hypotheses. Under the null, the asymptotic distribution of the maximum-type statistic depends on the spectra of the underlying kernels, with their effective ranks governing the transition between non-Gaussian chaos and Gaussian behavior.

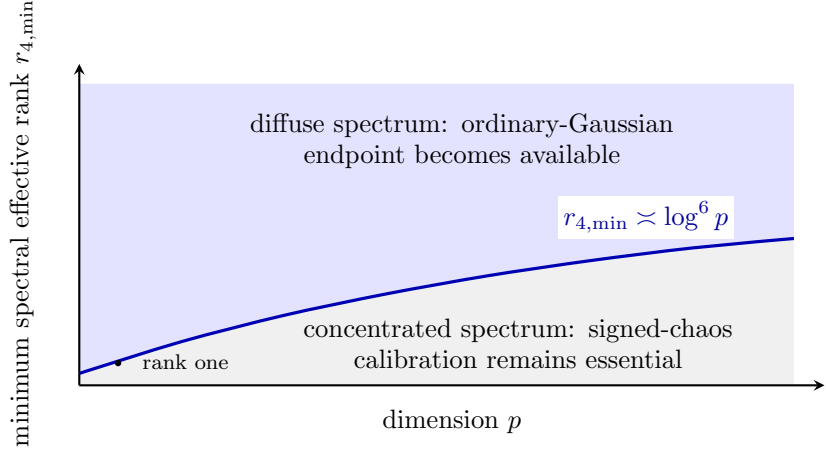
\begin{figure}[H]
\centering
\begingroup
\setstretch{1}
\begin{tikzpicture}[x=1.35cm,y=1.05cm,>=stealth,font=\small]
  \fill[gray!12] (0.8,0.7) rectangle (7.8,4.5);
  \fill[blue!11]
    (0.8,4.5)--(7.8,4.5)--(7.8,2.55)--(6.8,2.43)--(5.8,2.27)--
    (4.8,2.08)--(3.8,1.85)--(2.8,1.58)--(1.8,1.25)--(0.8,0.85)--cycle;
  \draw[->,thick] (0.8,0.7)--(8.1,0.7);
  \node at (4.45,0.25) {dimension $p$};
  \draw[->,thick] (0.8,0.7)--(0.8,4.75);
  \node[rotate=90] at (0.25,2.75)
    {minimum spectral effective rank $r_{4,\min}$};
  \draw[blue!70!black,very thick]
    plot[smooth] coordinates {(0.8,0.85) (1.8,1.25) (2.8,1.58)
      (3.8,1.85) (4.8,2.08) (5.8,2.27) (6.8,2.43) (7.8,2.55)};
  \node[blue!60!black,fill=white,inner sep=2pt] at (6.35,2.82)
    {$r_{4,\min}\asymp\log^6 p$};
  \node[align=center,text width=6cm] at (4.55,3.78)
    {diffuse spectrum: ordinary-Gaussian\\endpoint becomes available};
  \node[align=center,text width=6cm] at (5.1,1.23)
    {concentrated spectrum: signed-chaos\\calibration remains essential};
  \fill (1.18,0.98) circle (1.2pt);
  \node[anchor=west,font=\scriptsize] at (1.32,1.00) {rank one};
\end{tikzpicture}
\endgroup
\caption{Schematic phase transition indexed by dimension and spectral
diffuseness.  The fourth-order effective rank is built from the eigenvalues
through $r_{4,j}=\left(\sum_r\lambda_{jr}^2\right)^2/\sum_r\lambda_{jr}^4$.  Bounded
$r_{4,\min}$ retains a non-Gaussian signed-chaos target, whereas
$r_{4,\min}/\log^6p\to\infty$ is sufficient for the ordinary-Gaussian
endpoint once the finite-sample comparison terms also vanish.  The curve is
schematic and not to scale.}
\label{fig:phase-transition}
\end{figure}

The main contributions are summarized as follows.
\begin{enumerate}
\renewcommand{\labelenumi}{\textup{(\roman{enumi})}}
\item We construct a joint signed-chaos target and establish a L\'evy approximation result for maxima of canonical order-two \(U\)-statistics, allowing arbitrary signed spectra, unrelated eigensystems, and finite or infinite rank, without introducing a spectral truncation remainder. The key point is that the approximation preserves the full cross-coordinate dependence even when the underlying kernels do not admit a common eigenbasis.

\item We derive sharper anti-concentration bounds by exploiting additional spectral structure, which further yields explicit distributional approximation rates for maxima of canonical  \(U\)-statistics in high dimensions, allowing the dimension \(p\) to grow exponentially with the sample size.

\item We identify a   phase transition driven by the effective rank from non-Gaussian signed-chaos behavior to an ordinary Gaussian maximum; see Figure~\ref{fig:phase-transition}. Specifically, concentrated spectra retain the non-Gaussian   limit, whereas sufficiently diffuse spectra lead to Gaussianization toward the covariance-matched Gaussian maximum.

\item For positive-semidefinite kernels, we propose and establish the validity of a Gaussian multiplier bootstrap for practical implementation. The procedure requires neither estimation of  eigensystems, nor prior knowledge of whether the relevant approximation is a non-Gaussian chaos or an ordinary Gaussian law.

\end{enumerate}

The rest of the paper is organized as follows.
Section~\ref{sec:setup} formulates the problem and develops the approximation theory.  Section~\ref{sec:phase-transition} gives the spectral
phase transition.  Section~\ref{sec:feasible-calibration} introduces the proposed
bootstrap and proves its validity, and Section~\ref{sec:applications} gives the two
statistical applications.  Section~\ref{sec:simulations} reports the numerical
study.  Technical proofs and additional numerical results are provided in the
supplementary material.

\section{Gaussian chaos approximation}\label{sec:setup}

\subsection{Preliminaries}\label{sec:pre}
Throughout this paper, for hyperrectangles $\cR_p$ in $\R^p$ and random vectors $\bm Y,\bm Z\in\mathbb R^p$, define
\begin{equation*}
 \rho_{\cR}(\bm Y,\bm Z)
 =\sup_{A\in\cR_p}\left|\Pp(\bm Y\in A)-\Pp(\bm Z\in A)\right|.
\end{equation*}
For  random variables $\xi$ and $\zeta$, write
\begin{equation*}
 \dK(\xi,\zeta)=\sup_{t\in\R}|\Pp(\xi\le t)-\Pp(\zeta\le t)|,
\end{equation*}
\begin{equation*}
 \dL(\xi,\zeta)=\inf\{u>0:\Pp(\zeta\le t-u)-u\le\Pp(\xi\le t)\le\Pp(\zeta\le t+u)+u\ \text{for all }t\}.
\end{equation*}
The concentration modulus of $\xi$ is defined by
\begin{equation*}
 \omega_{\xi}(u)=\sup_{t\in\R}\Pp(t<\xi\le t+u),\qquad u>0.
\end{equation*}

For a symmetric canonical kernel $h$, let
$\mathcal T_h:\mathcal L^2(\mathbb P)\to\mathcal L^2(\mathbb P)$ be the Hilbert--Schmidt operator
\begin{equation*}
 (\mathcal T_hf)(x)=\int h(x,y)f(y)\mathbb P(\dd y).
\end{equation*}
Its nonzero spectral expansion is
\begin{equation*}
 h(x,y)=\sum_{r=1}^{\infty}\lambda_r\phi_r(x)\phi_r(y)
 \quad\text{in }\mathcal L^2(\mathbb P^2),\qquad
 \|h\|_{\mathcal L^2(\mathbb P^2)}^2=\sum_{r=1}^{\infty}\lambda_r^2,
\end{equation*}
where the eigenfunction satisfies that $\E\phi_r(X)=0$.
   The first contraction satisfies
\begin{equation*}
 h\starone h=\sum_{r=1}^{\infty}\lambda_r^2\phi_r\otimes\phi_r,
 \qquad
 \|h\starone h\|_{\mathcal L^2(\mathbb P^2)}^2
 =\sum_{r=1}^{\infty}\lambda_r^4.
\end{equation*}
Define the fourth-order effective rank by
\begin{equation*}
 r_4(h):=\frac{\left(\sum_{r=1}^{\infty}\lambda_r^2\right)^2}
 {\sum_{r=1}^{\infty}\lambda_r^4},
\end{equation*}
which
 equals $1$ for a rank-one kernel and equals $m$ when the spectrum has $m$
eigenvalues of equal magnitude.
A large effective rank indicates that the variance is spread over many eigendirections.
 Consider the normalized second-chaos variable
$ Q_h=\left\{2\sum_{r=1}^{\infty}\lambda_r^2\right\}^{-1/2}
\\\sum_{r=1}^{\infty}\lambda_r(G_r^2-1)$, where $G_1,G_2,\ldots$ are independent standard Gaussian variables. One has
$\E Q_h^2=1$ and
$\operatorname{cum}_4(Q_h):=\E Q_h^4-3=12/r_4(h).$
The identity
$\operatorname{cum}_4(Q_h)=12/r_4(h)$ further makes this geometry the precise
fourth-moment obstruction to Gaussianity, in the same spirit as classical
Gaussian-limit criteria for degenerate $U$-statistics and Gaussian chaoses
\cite{deJong1987,NourdinPeccatiReinert2010,DoblerPeccati2017,
DoblerPeccati2019,Koike2023}.

\subsection{Problem formulation}\label{sec:Problem settings}
For each $n\ge2$, let $X_1,\ldots,X_n$ be independent and identically
distributed random elements of a measurable space $(\cX,\cA)$ with common
law $\mathbb P$.  Let $p=p_n\ge3$, and let
$h_1,\ldots,h_p\in\mathcal L^2(\mathbb P^2)$ be symmetric and canonical,
so that $\E\left\{h_j(x,X)\right\}=0$ for $\mathbb P$-almost every $x$.
All objects may depend on $n$.  Assume $a_j:=\|h_j\|_{\mathcal L^2(\mathbb P^2)}>0$, and define
\begin{equation*}
 U_{n,j}=\frac{2}{n(n-1)}\sum_{1\le i<i'\le n}h_j(X_i,X_{i'}),
 \qquad
 W_{n,j}=\sqrt{\frac{n(n-1)}{2}}\frac{U_{n,j}}{a_j}.
\end{equation*}
  Put
\begin{equation*}
 \bm W_n=\left(W_{n,1},\ldots,W_{n,p}\right)^\top,
 \qquad
 M_n=\max_{j\le p}W_{n,j}.
\end{equation*}
Let $\bm\Sigma=\left(\Sigma_{jk}\right)_{j,k\le p}$, where
$\Sigma_{jk}:=\Cov(W_{n,j},W_{n,k})
 ={(a_ja_k)^{-1}}{\langle h_j,h_k\rangle_{\mathcal L^2(\mathbb P^2)}}$.
Here, $\mathbf\Sigma$ is allowed to be a singular covariance matrix.





For every $1\le j\le p$,
\begin{equation*}
 h_j(x,y)=\sum_{r=1}^{\infty}\lambda_{jr}\phi_{jr}(x)\phi_{jr}(y)
 \quad\text{in }\mathcal L^2(\mathbb P^2),
 \qquad
 a_j^2=\sum_{r=1}^{\infty}\lambda_{jr}^2,
\end{equation*}
where
$|\lambda_{j1}|\ge|\lambda_{j2}|\ge\cdots\ge0$. Here, the eigenvalues
$\lambda_{jr}$ may have arbitrary signs.  For each fixed $j$, the
eigenfunctions corresponding to nonzero eigenvalues are orthonormal.  Let
\(\cH_0=\overline{\operatorname{span}}
 \{\phi_{jr}:j\le p,\ r\ge1,\ \lambda_{jr}\ne0\}
 \subset \mathcal L_0^2(\mathbb P)\)
and let $\mathbb G_0(\cdot)$ be an isonormal Gaussian process over $\cH_0$. Let $Z_{jr}=\mathbb G_0(\phi_{jr})$ and  $\E(Z_{jr}Z_{ks})
 =\langle\phi_{jr},\phi_{ks}\rangle_{\mathcal L^2(\mathbb P)}$, and define the joint signed-chaos target and its maximum by
\begin{equation*}
 Q_j=\frac1{\sqrt2a_j}\sum_{r=1}^{\infty}
 \lambda_{jr}(Z_{jr}^2-1),
 \qquad
 M_Q=\max_{j\le p}Q_j.
\end{equation*}
The above series converges in $\mathcal L^2$. Besides,
\begin{equation*}
 \Cov(Q_j,Q_k)=\Cov(W_{n,j},W_{n,k})
 =\frac{\langle h_j,h_k\rangle_{\mathcal L^2(\mathbb P^2)}}{a_ja_k}
 =\frac1{a_ja_k}\sum_{r=1}^{\infty}\sum_{s=1}^{\infty}
 \lambda_{jr}\lambda_{ks}
 \langle\phi_{jr},\phi_{ks}\rangle_{\mathcal L^2(\mathbb P)}^2,
\end{equation*}
in which the double series is absolutely convergent.  Moreover, every $Q_j$ and hence
$M_Q$ has an atomless distribution.

The use of a common isonormal Gaussian process $\mathbb G_0(\cdot)$ is essential for preserving the joint dependence across coordinates. If each kernel were diagonalized separately and equipped with an independent collection of Gaussian scores, the resulting construction would reproduce the marginal distribution of each \(Q_j\) but would generally fail to preserve the cross-coordinate covariance structure. By contrast, the common isonormal process ensures that
\(
\mathbb E(Z_{jr}Z_{ks})
=
\langle \phi_{jr},\phi_{ks}\rangle_{\mathcal L^2(P)}\),
so that the kernels need not share a common eigensystem, and arbitrary  cross-coordinate dependence can be accommodated without simultaneous diagonalization.

\subsection{Approximation theory}\label{sec:all-rank}

For a random variable $Y$ and $\nu\in\{1,2\}$, define $
 \|Y\|_{\psi_\nu}
 =\allowbreak\inf\{c>0:\E\exp(|Y|^\nu/c^\nu)\le2\}.
$ Let $X'$ denote an independent copy of $X$ and define the leverage function
\begin{equation*}
 \ell_j(x)=\left\{\E\left[\left\{h_j(x,X')/a_j\right\}^2\right]\right\}^{1/2}
 =\left\|h_j(x,\cdot)/a_j\right\|_{\mathcal L^2(\mathbb P)}.
\end{equation*}
Then,  one has $\E\ell_j(X)^2=1$.
{\begin{assumption}\label{ass:kernel-leverage}
Fix exponents $\alpha_h,\alpha_\ell\ge0$.  There are finite constants
\begin{equation*}
 K_h=\max_{j\le p}\sup_{{q\in\N,q\ge2}}
       q^{-\alpha_h}\|h_j(X,X')/a_j\|_{\mathcal L^q},
 \qquad
 K_\ell=\max_{j\le p}\sup_{{q\in\N,q\ge2}}
       q^{-\alpha_\ell}\|\ell_j(X)\|_{\mathcal L^q}.
\end{equation*}
In the following, denote $
 \mathcal K=1+K_h+K_\ell,
 \varrho=\max\{1,\alpha_\ell+1/2\},
 \varkappa=\max\{2,\alpha_h+1,\alpha_\ell+3/2\}.$

\end{assumption}
}

\begin{remark}\label{rem:no-high-moment}
The quantity
$K_h$ bounds an individual standardized kernel value, while $K_\ell$
controls how strongly one observation can influence the quadratic statistic
after its partner is integrated out.  The common
sub-exponential/sub-Gaussian regime corresponds to
$(\alpha_h,\alpha_\ell)=(1,1/2)$, and uniformly bounded standardized kernels
$h_j/a_j$ satisfy the assumption with $(0,0)$.  Products of Gaussian
covariance scores satisfy $(2,1)$; see
Section~\ref{sec:gaussian-covariance-app}.   Similar moment conditions are also imposed in the literature of Gaussian
approximations for maxima of sums \cite{CCK2013,CCK2017,CCKK2022}.
\end{remark}


The distributional approximation problem has two  parts.  The first is a
comparison between the maximum of degenerate U-statistics $M_n$ and the maximum of the joint signed-chaos target $M_Q$.  The
second controls the concentration of $M_Q$ over short intervals. Generally, the
L\'evy result below needs only the first part, whereas the Kolmogorov conclusion
also needs the second.

We first provide the L\'evy approximation results.

For $0<\eta<1$, let $
 r_{n,p}(\eta)
 = \sqrt{\log(p/\eta)/n}
 +\log(p/\eta)/n,
$
and define
\begin{equation}
 \Delta_n
 = \mathcal K^3\log^2(p)
 \left\{\frac{\log^{3\varrho}(np)}{\sqrt n}
       +\frac{\log^{3\varkappa}(np)}{n^2}\right\}.
 \label{eq:Delta-n-signed}
\end{equation}
\begin{proposition}
\label{thm:all-rank-chaos}
Suppose Assumption~\ref{ass:kernel-leverage} holds.  For every $0<\eta<1$,
\begin{equation*}
 \dL(M_n,M_Q)
 \lesssim \Delta_n^{1/4}+r_{n,p}(\eta)+\eta.
\end{equation*}
\end{proposition}

The proposition  is general and applies to kernels of finite or infinite rank with arbitrary signed spectra. It is also the basic theoretical bridge: it replaces the
data maximum $M_n$ by the maximum of the joint chaos $M_Q$ without asking whether that chaos is
 close to a Gaussian vector or not.
 L\'evy
distance tolerates a small horizontal shift and therefore needs no
anti-concentration assumption, so the result remains valid even when a
shrinking interval   carries substantial target probability.  In particular, from
\eqref{eq:Delta-n-signed},
\begin{equation*}
 \dL(M_n,M_Q)
 \lesssim
 \mathcal K^{3/4}
 \left\{
 \frac{\log^{(3\varrho+2)/4}(np)}{n^{1/8}}
 +\frac{\log^{(3\varkappa+2)/4}(np)}{n^{1/2}}
 \right\}
 +r_{n,p}(\eta)+\eta.
\end{equation*}
Hence, for bounded $\mathcal K$, the baseline regimes
$(\alpha_h,\alpha_\ell)=(1,1/2)$ and $(0,0)$ permit
$\log p=o(n^{1/10})$, while the Gaussian-product regime
$(\alpha_h,\alpha_\ell)=(2,1)$ permits $\log p=o(n^{1/13})$.

\begin{proposition}
\label{prop:Levy-Kolmogorov-transfer}
Let $\xi$ and $\zeta$ be real random variables.  For every
$\delta>\dL(\xi,\zeta)$,
\begin{equation*}
 \dK(\xi,\zeta)\le \delta+\omega_\zeta(2\delta).
\end{equation*}
\end{proposition}
Proposition \ref{prop:Levy-Kolmogorov-transfer} shows that a L\'evy bound can be upgraded to a Kolmogorov bound whenever the approximating law has sufficient anti-concentration at the corresponding error scale.
A Kolmogorov bound can be obtained by directly combining the L\'evy approximation in Proposition \ref{thm:all-rank-chaos} with Proposition \ref{prop:Levy-Kolmogorov-transfer}. However, a sharper result is available by retaining the smoothing parameter \(\varepsilon\) in the comparison step and optimizing it jointly with the anti-concentration bound for \(M_Q\), as shown in the following theorem.
\begin{theorem}
\label{thm:all-rank-Kolmogorov}
Suppose Assumption~\ref{ass:kernel-leverage} holds.  For every $0<\eta<1$ and
$\eps>0$,
\begin{equation}
 \dK(M_n,M_Q)
 \lesssim \Delta_n\eps^{-3}
 +\omega_{M_Q}\left(C\{\eps+r_{n,p}(\eta)\}\right)
 +\eta.
 \label{eq:all-rank-Kolmogorov}
\end{equation}
\end{theorem}

The first term on the right-hand side of \eqref{eq:all-rank-Kolmogorov} measures the stochastic replacement error, while the concentration modulus quantifies the sensitivity of the target distribution to perturbations of the critical value induced by the approximation errors. Consequently, any available anti-concentration bound for \(M_Q\) can be inserted into \eqref{eq:all-rank-Kolmogorov} to obtain an explicit bound on the Kolmogorov approximation error.


\subsection{Anti-concentration bounds}

The following proposition shows that, in the absence of additional structural assumptions, the anti-concentration bound necessarily deteriorates linearly with the dimension \(p\).


\begin{proposition}[General anti-concentration and sharpness]
\label{prop:general-chaos-anticoncentration}
Let $Q_1,\ldots,Q_p$ be jointly defined centered, variance-one elements of the
second homogeneous Gaussian chaos, with arbitrary dependence, signs, and
ranks.  Then
\begin{equation}
 \omega_{M_Q}(u)\lesssim p\sqrt u,
 \qquad u>0.
 \label{eq:general-chaos-anti}
\end{equation}
Conversely, for independent standard Gaussian variables
$Z_1,\ldots,Z_p$ and $Q_j=(1-Z_j^2)/\sqrt2$, there are universal constants
$c,u_0>0$ such that
\begin{equation*}
 \omega_{M_Q}(u)\ge c\{1\wedge p\sqrt u\},
 \qquad 0<u\le u_0.
\end{equation*}
Consequently, $1\wedge p\sqrt u$ is the sharp uniform order over this class.
\end{proposition}

The matching lower bound is attained by an explicit negative rank-one construction, in which the maximum can  accumulate substantial probability near the common upper endpoints of the negative rank-one coordinates.
 Thus, the dependence on \(p\) cannot be improved without imposing further spectral structure.


We therefore turn to positive-spectrum settings with additional spectral structure, including finite-rank kernels and infinite-rank kernels with polynomial or geometric eigenvalue decay, for which  sharper anti-concentration bounds can be obtained.

\begin{assumption}
\label{ass:block-PSD}
For every $j\le p$, there are a real separable Hilbert space $\mathcal H_j$ and a
strongly measurable feature map $\Phi_j:\cX\to\mathcal H_j$ such that
\begin{equation*}
 h_j(x,y)=\langle\Phi_j(x),\Phi_j(y)\rangle_{\mathcal H_j},
 \qquad
 \E\Phi_j(X)=0,
 \qquad
 \E\|\Phi_j(X)\|_{\mathcal H_j}^2<\infty.
\end{equation*}
For $j,k\le p$, define the normalized block covariance operator
\begin{equation*}
 \mathcal C_{jk}
 =\frac{\E\left\{\Phi_j(X)\otimes\Phi_k(X)\right\}}
 {\sqrt{a_ja_k}}:
 \mathcal H_k\to\mathcal H_j.
\end{equation*}
Then $\|\mathcal C_{jj}\|_{\mathrm{HS}}=1$.  Let
$\bm G=\left(G_1,\ldots,G_p\right)$ be one centered jointly Gaussian block
vector with covariance blocks $\mathcal C_{jk}$, and let $
 \tau_j=\tr(\mathcal C_{jj}),
 \tau_- = \min_{j\le p}\tau_j,
 \tau_+ = \max_{j\le p}\tau_j. $
If $\mu_{j1}\ge\mu_{j2}\ge\cdots\ge0$ are the eigenvalues of
$\mathcal C_{jj}$, then $\sum_{r=1}^{\infty}\mu_{jr}^2=1$,  $\tau_j=\sum_{r=1}^{\infty}\mu_{jr}$,
and the positive-spectrum chaos has the joint representation
\begin{equation}
 Q_j\overset d=\frac{\|G_j\|^2-\tau_j}{\sqrt2}
 =\frac1{\sqrt2}\sum_{r=1}^{\infty}\mu_{jr}\left(Z_{jr}^2-1\right),
 \label{eq:PSD-Gaussian-block-representation}
\end{equation}
where $\{Z_{jr}\}_{r=1}^{\infty}$ are independent standard Gaussian variables.
\end{assumption}

\begin{theorem}
\label{thm:structured-anticoncentration}
Suppose that Assumption \ref{ass:block-PSD} holds.
\begin{enumerate}
\renewcommand{\labelenumi}{\textup{(\roman{enumi})}}
\item Suppose that $\rank(\mathcal C_{jj})\le M$ for every $j\le p$.
Equivalently, after ordering and zero padding, $\mu_{jr}=0$ for $r>M$.
Then, for every $u>0$,
\begin{equation}
 \omega_{M_Q}(u)
 \lesssim
 \left\{\tau_- M u
 \log\left(\frac{ep\sqrt{\tau_+}}{u\wedge1}\right)\right\}^{1/3}
 +\sqrt{\tau_-u}.
 \label{eq:finite-rank-PSD-anti}
\end{equation}

\item Suppose   that every $\mathcal C_{jj}$ has infinite rank and, uniformly
in $j$ and $r$,
\[
 c_2r^{-2\beta}\le\mu_{jr}\le c_1r^{-2\alpha},
 \qquad \beta\ge\alpha>1/2.
\]  Then, for
$0<\rho\le1$ and $0<u\le1/2$,
\begin{equation*}
 \begin{aligned}
 \sup_{\{t:\,\tau_-+\sqrt2t\ge\rho\}}
 \Pp(t<M_Q\le t+u)
  &\lesssim \rho^{-1/2-2\beta/(2\alpha-1)}u
  \log^{1/2+2\beta/(2\alpha-1)}(ep)\\
  &\qquad\times\log^{2\beta/(2\alpha-1)}(e/u),
 \end{aligned}
\end{equation*}
where the implicit constant depends only on $\alpha,\beta,c_1$, and $c_2$.

\item Suppose that every $\mathcal C_{jj}$ has infinite rank and, uniformly
in $j$ and $r$,
\[
 c_2e^{-\beta r}\le\mu_{jr}\le c_1e^{-\alpha r},
 \qquad \beta\ge\alpha>0.
\]
Then, for $0<\rho\le1$ and $0<u\le1/2$,
\begin{equation*}
 \begin{aligned}
 \sup_{\{t:\,\tau_-+\sqrt2t\ge\rho\}}
 \Pp(t<M_Q\le t+u)
 &\lesssim \rho^{-1/2-\beta/\alpha}u
 \log^{1/2+\beta/\alpha}(ep)
 \log^{\beta/\alpha}(e/u),
 \end{aligned}
\end{equation*}
where the implicit constant depends only on $\alpha,\beta,c_1$, and $c_2$.
\end{enumerate}
\end{theorem}

The finite-rank result improves the linear dependence on \(p\) in the general bound to logarithmic dependence. The key gain comes from positivity: each coordinate can be represented as a centered squared Gaussian norm, while finite rank restricts the number of Gaussian directions contributing to its fluctuations. The heterogeneous lower endpoints are accommodated through \(\tau_-\) and \(\tau_+\) in \eqref{eq:finite-rank-PSD-anti}. For infinite-rank kernels, analogous improvements remain available under regular spectral decay. The lower spectral envelope ensures sufficiently many nonnegligible directions, whereas the upper envelope controls the contribution of the spectral tail. The resulting bounds are local away from the lower edge \(-\tau_-/\sqrt{2}\), near which the density of a centered squared Gaussian norm may become singular.


\begin{corollary}\label{cor:explicit-all-rank}
Suppose Assumptions~\ref{ass:kernel-leverage} and \ref{ass:block-PSD} hold.
\begin{enumerate}
\renewcommand{\labelenumi}{\textup{(\roman{enumi})}}
\item Under Theorem~\ref{thm:structured-anticoncentration}(i),
\begin{equation*}
 \begin{aligned}
 \dK(M_n,M_Q)
 &\lesssim
  \Delta_n^{1/10}\left\{\tau_-M\Lambda_{n,p}(\eta)\right\}^{3/10}
  +\left\{\tau_-M r_{n,p}(\eta)\Lambda_{n,p}(\eta)\right\}^{1/3}\\
  &\quad+\sqrt{\tau_-r_{n,p}(\eta)}+\eta,
 \end{aligned}
\end{equation*}
where \begin{equation*}
 \Lambda_{n,p}(\eta)
 =\log\left\{
 \frac{ep\sqrt{\tau_+}}
 {\left\{\left(\Delta_n+r_{n,p}(\eta)\right)\wedge1\right\}}
 \right\}.
\end{equation*}

\item Let \begin{equation*}
 B_{n,\rho}^{\rm pol}(\eta)
 =\rho^{-1/2-2\beta/(2\alpha-1)}
  \log^{1/2+2\beta/(2\alpha-1)}(ep)
  \log^{2\beta/(2\alpha-1)}\!\left(
  \frac{e}{\left\{\left(\Delta_n+r_{n,p}(\eta)\right)\wedge1\right\}}
  \right).
\end{equation*}
Under Theorem~\ref{thm:structured-anticoncentration}(ii), if $C\left[\left\{{\Delta_n}/{B_{n,\rho}^{\rm pol}(\eta)}\right\}^{1/4}
 +r_{n,p}(\eta)\right]\le\rho$ for some $0<\rho\le1$,  then
\begin{equation*}
 \begin{aligned}
 &\sup_{\left\{t:\,\tau_-+\sqrt2t\ge2\rho\right\}}
 \left|\Pp(M_n\le t)-\Pp(M_Q\le t)\right|\lesssim
 \Delta_n^{1/4}\{B_{n,\rho}^{\rm pol}(\eta)\}^{3/4}
 +r_{n,p}(\eta)B_{n,\rho}^{\rm pol}(\eta)+\eta.
 \end{aligned}
\end{equation*}

\item Let
\begin{equation*}
 B_{n,\rho}^{\rm geo}(\eta)
 =\rho^{-1/2-\beta/\alpha}
  \log^{1/2+\beta/\alpha}(ep)
  \log^{\beta/\alpha}\!\left(
  \frac{e}{\left\{\left(\Delta_n+r_{n,p}(\eta)\right)\wedge1\right\}}
  \right).
\end{equation*}
Under Theorem~\ref{thm:structured-anticoncentration}(iii),
if
$
 C\left[\left\{ {\Delta_n}/{B_{n,\rho}^{\rm geo}(\eta)}\right\}^{1/4}
 +r_{n,p}(\eta)\right]\le\rho
$
for some $0<\rho\le1$,
then
\begin{equation*}
 \begin{aligned}
 &\sup_{\left\{t:\,\tau_-+\sqrt2t\ge2\rho\right\}}
 \left|\Pp(M_n\le t)-\Pp(M_Q\le t)\right|\lesssim
 \Delta_n^{1/4}\{B_{n,\rho}^{\rm geo}(\eta)\}^{3/4}
 +r_{n,p}(\eta)B_{n,\rho}^{\rm geo}(\eta)+\eta.
 \end{aligned}
\end{equation*}
\end{enumerate}
\end{corollary}

Combining Theorem \ref{thm:all-rank-Kolmogorov} with the sharper anti-concentration bounds above and optimizing the smoothing parameter yields Corollary \ref{cor:explicit-all-rank}.
 The resulting rates differ across the spectral regimes through the corresponding anti-concentration bounds available under each set of assumptions.

\section{Spectral phase transition}\label{sec:phase-transition}

We now allow signed spectra and ask when the chaos itself admits an ordinary
Gaussian approximation.

For every $j\le p$, define
\begin{equation*}
 r_{4,j}
 =\frac{\left(\sum_{r=1}^{\infty}\lambda_{jr}^2\right)^2}{\sum_{r=1}^{\infty}\lambda_{jr}^4},
 \qquad
 r_{4,\min}=\min_{j\le p}r_{4,j}.
\end{equation*}
Let
$\bm Z_n=(Z_{n,1},\ldots,Z_{n,p})^\top\sim N(0,\bm\Sigma)$, where
$\Sigma_{jk}=\langle h_j,h_k\rangle_{\mathcal L^2(\mathbb P^2)}\allowbreak/(a_ja_k)$ is defined in Section~\ref{sec:Problem settings}. Here,
the covariance matrix $\mathbf\Sigma$ may be singular.

\begin{theorem}
\label{thm:chaos-Gaussianization}
For the general signed-chaos target,
\begin{equation*}
 \dK\left(M_Q,\max_{j\le p}Z_{n,j}\right)
 \lesssim   r_{4,\min}^{-1/6}\log(p),
\end{equation*}
and hence, for every $u>0$,
\begin{equation*}
 \omega_{M_Q}(u)
 \lesssim u\sqrt{\log(p)}
 +r_{4,\min}^{-1/6}\log(p).
\end{equation*}
\end{theorem}
The result permits arbitrary spectral signs and infinite rank.
 It implies that, if \(r_{4,\min}/\allowbreak\log^6 p\to\infty\), then the maximum of the signed-chaos vector can be further approximated by the maximum of a \(p\)-dimensional Gaussian vector with the same covariance structure, yielding the spectral transition summarized in Figure~\ref{fig:phase-transition}.  Since the fourth cumulant of coordinate
$j$ equals $12/r_{4,j}$ as discussed in Section \ref{sec:pre}, the worst coordinate controls the simultaneous
Gaussianization error.  The logarithmic factor is the price of comparing the
maximum of $p$ coordinates.


\begin{theorem}
\label{thm:two-stage-phase}
Suppose Assumption~\ref{ass:kernel-leverage} holds.  Combining
Theorems~\ref{thm:all-rank-Kolmogorov} and~\ref{thm:chaos-Gaussianization} and taking
$\eta=n^{-1}$ gives
\begin{align}
 \dK\left(M_n,\max_{j\le p}Z_{n,j}\right)
 \lesssim  {}
 &\mathcal K^{3/4}
 \left\{\frac{\log^{(6\varrho+7)/8}(np)}{n^{1/8}}
       +\frac{\log^{(6\varkappa+7)/8}(np)}{n^{1/2}}\right\}
 \notag\\
 &+\frac{\log(np)}{\sqrt n}
 +\frac{\log^{3/2}(np)}n
 +{r_{4,\min}^{-1/6}\log p}.
 \label{eq:two-stage-phase-bound}
\end{align}

\end{theorem}
We therefore summarize the phase transition as follows:
\begin{enumerate}
\renewcommand{\labelenumi}{\textup{(\roman{enumi})}}
\item If \(r_{4,\min}\) remains bounded, an ordinary Gaussian approximation is not valid uniformly over the model class, and the signed quadratic-chaos law remains the appropriate target. Valid critical-value approximation then requires sufficient anti-concentration of the chaos maximum at the relevant approximation scale, as quantified by Theorem~\ref{thm:all-rank-Kolmogorov}.

\item \begingroup\emergencystretch=1em
If \(r_{4,\min}/\log^6 p\to\infty\),
the maximum \(M_Q\) admits a Gaussian approximation by \(\max_{j\le p}Z_{n,j}\). When the finite-sample approximation terms also vanish, this yields an ordinary Gaussian approximation for the data maximum, as quantified by Theorem~\ref{thm:two-stage-phase}.
\par\endgroup
\end{enumerate}

\begin{remark}
In the common sub-exponential-kernel/sub-Gaussian-leverage regime,
$(\varrho,\varkappa)=(1,2)$, $\log p=o(n^{1/13})$ is sufficient for the
finite-sample part in Theorem \ref{thm:two-stage-phase} to vanish.
A bounded-kernel specialization of \citet[Theorem~2]{ImaiKoike2025} gives
the following benchmark:
\[
\dK\left(M_n,\max_{j\le p}Z_{n,j}\right)
 \lesssim
 r_{4,\min}^{-1/4}(\log p)^{3/2}
 +\left\{\frac{B_n^4\log^7(enp)}n\right\}^{1/4},
\]
where $B_n=\max_j\|h_j/a_j\|_\infty$.
Accordingly, for bounded kernels, the two results require \(\log p=o(n^{1/13})\) and \(\log p=o(n^{1/7})\), respectively, while both use the sufficient effective-rank condition \(r_{4,\min}/(\log p)^6\to\infty\). The slower rate of our result reflects its more general two-stage construction, which first approximates the statistic by the signed-chaos target and then Gaussianizes that target, with additional losses from smoothing and anti-concentration,  while \cite{ImaiKoike2025} bypasses the intermediate chaos approximation and directly approximates the maximum of the degenerate U-statistics by a Gaussian maximum.
\end{remark}

\section{Gaussian multiplier bootstrap}\label{sec:feasible-calibration}

\subsection{Construction and algorithm}
\label{sec:gram-method}
\begingroup\emergencystretch=1em
The  proposed Gaussian multiplier bootstrap approximates   the joint chaos distribution without
estimating eigensystems.
Let
$\mathbf H_j=(h_j(X_i,X_{i'}))_{i,i'\le n}$ be the symmetric Gram matrix, $\mathbf P_n=\mathbf I_n-n^{-1}\mathbf 1\mathbf 1^\top$, and
\begin{equation}
 \mathbf B_j=\mathbf P_n\mathbf H_j\mathbf P_n.
 \label{eq:gram-centered-matrix}
\end{equation}
\par\endgroup
Draw $\bm\xi^{(b)}\sim N(0,\mathbf I_n)$ independently of the data (use the same draw for
all coordinates), and compute
\begin{equation}
 Q_{n,j}^{\#,(b)}=
 \frac{(\bm\xi^{(b)})^\top \mathbf B_j\bm\xi^{(b)}-\tr(\mathbf B_j)}
 {\sqrt2\,\|\mathbf B_j\|_{\mathrm F}},
 \qquad T_n^{\#,(b)}=\max_{j\le p}Q_{n,j}^{\#,(b)}.
 \label{eq:gram-multiplier}
\end{equation}
Conditionally on the sample, each coordinate has mean zero and variance one.
The common multiplier vector preserves the empirical dependence among
kernels.

A scale estimator of $a_j$ is
\begin{equation*}
 \widehat a_{j,\mathrm{pair}}^2
 =\frac1{\lfloor n/2\rfloor}\sum_{r=1}^{\lfloor n/2\rfloor}
 h_j(X_{2r-1},X_{2r})^2,
 \qquad
 \widehat W_{n,j}
 =\sqrt{\frac{n(n-1)}{2}}
 \frac{U_{n,j}}{\widehat a_{j,\mathrm{pair}}},
\end{equation*}
where the disjoint pair terms are independent across $r$.  We note that the statistic
$U_{n,j}$ still uses every unordered pair; only its scale is estimated
from disjoint pairs.

\begin{algorithm}[H]
\caption{Proposed centered-Gram multiplier calibration}
\label{alg:gram-bootstrap}
\begin{algorithmic}[1]
\Require Data $X_1,\ldots,X_n$; kernels $h_1,\ldots,h_p$; multiplier draws
$B$; level $\alpha$.
\State Compute $U_{n,j}$ and $\widehat a_{j,\mathrm{pair}}$ for every
$j=1,\ldots,p$.
\State Form $\mathbf B_j$ according to \eqref{eq:gram-centered-matrix} for
every $j$.
\For{$b=1,\ldots,B$}
  \State Draw one common $\bm\xi^{(b)}\sim N(0,\mathbf I_n)$.
  \State Compute $T_n^{\#,(b)}$ according to \eqref{eq:gram-multiplier}.
\EndFor
\State Let $\widehat c_{1-\alpha}$ be the empirical $(1-\alpha)$ quantile of
$T_n^{\#,(1)},\ldots,T_n^{\#,(B)}$.
\State Reject if $\max_{j\le p}\widehat W_{n,j}>\widehat c_{1-\alpha}$.
\Statex \textit{Two-sided version:} augment the kernel collection with
$-h_1,\ldots,-h_p$ before applying the same steps.
\end{algorithmic}
\end{algorithm}


\subsection{Bootstrap validity}

We first control the estimated scales, then compare the conditional
bootstrap law with the population chaos.

\begin{theorem}
\label{thm:uniform-studentization}
Suppose Assumption~\ref{ass:kernel-leverage} holds.  For $0<\eta<1/2$, with probability
at least $1-\eta$,
\begin{equation*}
 \max_{j\le p}
 \left|\frac{\widehat a_{j,\mathrm{pair}}^2}{a_j^2}-1\right|
 \lesssim_{\alpha_h} K_h^2\left\{
 \sqrt{\frac{\log(8np/\eta)}n}
 +\frac{\log^{2\alpha_h+1}(8np/\eta)}n\right\}
 =:s_n(\eta).
\end{equation*}
{Let $\mathcal E_s(\eta)$ denote the event in the preceding
display.  If $s_n(\eta)\le1/2$, then on $\mathcal E_s(\eta)$,}
\begin{equation*}
 \max_{j\le p}\left|\widehat W_{n,j}-W_{n,j}\right|
 \le 2s_n(\eta)\max_{j\le p}|W_{n,j}|.
\end{equation*}
\end{theorem}
When the  approximation error vanishes,
$\max_j|W_{n,j}|=O_\Pp\{\log(np)\}$; studentization is asymptotically negligible if
$s_n(\eta_n)\log(np)\to0$ for a sequence $\eta_n\downarrow0$. An all-pairs estimator of $a_j$ may be more efficient, but its analysis would additionally require controlling the Hoeffding remainder of the squared-kernel U-statistic.


In the following, we suppose that Assumption~\ref{ass:block-PSD} holds and retain the
notation $\Phi_j$, $\mathcal C_{jk}$, $G_j$, and $\tau_j$ from
Assumption~\ref{ass:block-PSD}.  In particular,
$\|\mathcal C_{jj}\|_{\mathrm{HS}}=1$ and
$Q_j\overset d=\left(\|G_j\|^2-\tau_j\right)/\sqrt2$.
Let $
 \bar\Phi_j=\frac1n\sum_{i=1}^n\Phi_j(X_i), $
and define
\begin{equation*}
 \widetilde{\mathcal C}_{jk}
 =\frac{1}{n\sqrt{a_ja_k}}\sum_{i=1}^n
 \left\{\Phi_j(X_i)-\bar\Phi_j\right\}
 \otimes
 \left\{\Phi_k(X_i)-\bar\Phi_k\right\}.
\end{equation*}
Whenever
$\min_{j\le p}\|\widetilde{\mathcal C}_{jj}\|_{\mathrm{HS}}>0$, let
\begin{equation*}
 \widehat{\mathcal C}_{jk}
 =\frac{\widetilde{\mathcal C}_{jk}}
 {\left\{\|\widetilde{\mathcal C}_{jj}\|_{\mathrm{HS}}
          \|\widetilde{\mathcal C}_{kk}\|_{\mathrm{HS}}\right\}^{1/2}},
 \qquad
 \Delta=\max_{j,k\le p}
 \|\widehat{\mathcal C}_{jk}-\mathcal C_{jk}\|_{\mathrm{HS}}.
\end{equation*}

\begin{theorem}\label{thm:gram-bootstrap}
Suppose Assumption~\ref{ass:block-PSD} holds and $\min_{j\le p}\|\widetilde{\mathcal C}_{jj}\|_{\mathrm{HS}}>0$. For every $t\in\R$ and $u>0$,
\begin{align*}
 \Pp(M_Q\le t-2u)-C\Delta\frac{\log^2(ep)}{u^2}
 &\le \Pp^\#\!\left(\max_jQ_{n,j}^{\#}\le t\right)\\
 &\le \Pp(M_Q\le t+2u)+C\Delta\frac{\log^2(ep)}{u^2}.
\end{align*}
Consequently,
\begin{equation*}
 d_{\mathrm L}\!\left\{
 \operatorname{Law}^{\#}\!\left(\max_jQ_{n,j}^{\#}\right),\operatorname{Law}(M_Q)\right\}
 \lesssim \{\Delta\log^2(ep)\}^{1/3},
\end{equation*}
\begin{equation}
 \sup_t\left|
 \Pp^\#\!\left(\max_jQ_{n,j}^{\#}\le t\right)-\Pp(M_Q\le t)
 \right|
 \le C\Delta\frac{\log^2(ep)}{u^2}
 +\omega_{M_Q}(2u).
 \label{eq:gram-kolmogorov}
\end{equation}
\begin{enumerate}
\renewcommand{\labelenumi}{\textup{(\roman{enumi})}}
\item Under the conditions of
Theorem~\ref{thm:structured-anticoncentration}(i), let
\[
 \Lambda=\log\left\{\frac{ep(1+\tau_+M)}{\Delta\log^2(ep)}\right\}.
\]
If $0<\Delta\log^2(ep)\le e^{-1}$, optimizing \eqref{eq:gram-kolmogorov} with
\eqref{eq:finite-rank-PSD-anti} gives the sharper global bound
\begin{equation*}
 \dK\!\left\{
 \operatorname{Law}^{\#}\!\left(\max_jQ_{n,j}^{\#}\right),\operatorname{Law}(M_Q)
 \right\}
 \lesssim \{\Delta\log^2(ep)\}^{1/7}\{\tau_-M\Lambda\}^{2/7}.
\end{equation*}

\item Under the polynomial spectral conditions of
Theorem~\ref{thm:structured-anticoncentration}(ii), let
\[
 B_\rho^{\rm pol}
 =\rho^{-1/2-2\beta/(2\alpha-1)}
  \log^{1/2+2\beta/(2\alpha-1)}(ep)
  \log^{2\beta/(2\alpha-1)}\!\left(\frac{e}{\{\Delta\log^2(ep)\}\wedge1}\right).
\]
If $
 C\left\{\Delta\log^2(ep)/B_\rho^{\rm pol}\right\}^{1/3}\le\rho\leq 1,$
then
\begin{equation*}
 \begin{aligned}
 &\sup_{\{t:\,\tau_-+\sqrt2t\ge2\rho\}}
 \left|\Pp^\#\!\left(\max_jQ_{n,j}^{\#}\le t\right)-\Pp(M_Q\le t)\right|\lesssim \{\Delta\log^2(ep)\}^{1/3}(B_\rho^{\rm pol})^{2/3}.
 \end{aligned}
\end{equation*}

\item Under the geometric spectral conditions of
Theorem~\ref{thm:structured-anticoncentration}(iii), let
\[
 B_\rho^{\rm geo}
 =\rho^{-1/2-\beta/\alpha}
  \log^{1/2+\beta/\alpha}(ep)
  \log^{\beta/\alpha}\!\left(\frac{e}{\{\Delta\log^2(ep)\}\wedge1}\right).
\]
If
\(
 C\left\{{\Delta\log^2(ep)}/{B_\rho^{\rm geo}}\right\}^{1/3}\le\rho\leq 1,
\)
then
\begin{equation*}
 \begin{aligned}
 &\sup_{\{t:\,\tau_-+\sqrt2t\ge2\rho\}}
 \left|\Pp^\#\!\left(\max_jQ_{n,j}^{\#}\le t\right)-\Pp(M_Q\le t)\right|\lesssim \{\Delta\log^2(ep)\}^{1/3}(B_\rho^{\rm geo})^{2/3}.
 \end{aligned}
\end{equation*}
\end{enumerate}

\end{theorem}


The smoothed comparison error is controlled by the uniform block covariance discrepancy
\(\Delta\). The L\'evy bound requires no density control, whereas the Kolmogorov bounds are obtained by combining the comparison with the anti-concentration results in Theorem~\ref{thm:structured-anticoncentration}. The use of common multipliers preserves the full cross-coordinate dependence. Section~\ref{sec:applications} further verifies the required covariance and spectral conditions for the two applications.

In addition, we note that for a general \(\mathcal L^2(\mathbb P^2)\) kernel, its values on the diagonal are not determined by its \(\mathcal L^2(\mathbb P^2)\) equivalence class. Since the centered-Gram bootstrap uses the diagonal entries \(h_j(X_i,X_i)\), we require a pointwise-defined kernel representative with a specified diagonal.
 The established theory above extends fixed-dimensional bootstrap results for degenerate \(U\)- and \(V\)-statistics and wild multipliers \cite{ArconesGine1992,LeuchtNeumann2013,Chwialkowski2014} to simultaneous maxima of PSD kernels, without requiring estimation of the underlying eigensystems.

\section{Statistical applications}
\label{sec:applications}

We consider simultaneous testing over covariance-graph edges and over
coordinates or blocks in a two-sample problem.

\subsection{Gaussian covariance-graph testing}
\label{sec:gaussian-covariance-app}
Let $\widetilde{\bm X}_1,\ldots,\widetilde{\bm X}_{2n}$ be independent
$N_d(\bm\mu,\bm\Sigma)$ vectors.  To remove an unknown mean without estimating it,
form disjoint differences
\[
 \bm X_i=\frac{\widetilde{\bm X}_{2i}-\widetilde{\bm X}_{2i-1}}{\sqrt2},
 \qquad i=1,\ldots,n,
\]
so $\bm X_i\sim N_d(\bm 0,\bm\Sigma)$.
Let
$\mathcal E\subset\{(j,k):1\le j<k\le d\}$ be a candidate edge set and let
$p=|\mathcal E|$.  The null hypothesis is
\[
 H_{0,\mathcal E}:\quad \sigma_{jk}=0
 \quad\text{for every }(j,k)\in\mathcal E.
\]
{
For each $(j,k)\in\mathcal E$, define
\begin{equation}
 h_{jk}(\bm x,\bm y)=x_jx_k\,y_jy_k,
 \qquad
 W_{n,jk}=
 \sqrt{\frac{2}{n(n-1)}}
 \frac{\sum_{i<\ell}X_{ij}X_{ik}X_{\ell j}X_{\ell k}}
 {\E(X_j^2X_k^2)}.
 \label{eq:covariance-U-kernel}
\end{equation}
The pair average estimates the squared edge signal because
$\E h_{jk}(\bm X,\bm X')=\sigma_{jk}^2$.  Thus the procedure is a quadratic analogue of
maximum sample-correlation and coherence tests used for sparse covariance
structure \cite{CaiJiang2011,CaiLiuXia2013}.

\begin{proposition}
\label{prop:gaussian-covariance-reduction}
Under $H_{0,\mathcal E}$, every kernel in
\eqref{eq:covariance-U-kernel} is canonical and has rank one.  With
$\psi_{jk}(\bm x)=x_jx_k/\{\E(X_j^2X_k^2)\}^{1/2}$,
${h_{jk}(\bm x,\bm y)}/{\E(X_j^2X_k^2)}
 =\psi_{jk}(\bm x)\psi_{jk}(\bm y),$ and $ Q_{jk}= {(Z_{jk}^2-1)}/{\sqrt2}$,
where $\bm Z=\left(Z_{jk}\right)_{(j,k)\in\mathcal E}$ is Gaussian with
\[
\E(Z_{jk}Z_{j'k'})=\E\{\psi_{jk}(\bm X)\psi_{j'k'}(\bm X)\}.
\]
Consequently,
\begin{equation*}
 \max_{(j,k)\in\mathcal E}Q_{jk}
 =\frac{\max_{(j,k)\in\mathcal E}|Z_{jk}|^2-1}{\sqrt2},
 \qquad r_{4,jk}=1.
\end{equation*}
Uniformly over the edge set,
$\|h_{jk}(X,X')/\E(X_j^2X_k^2)\|_{\mathcal L^q}\le Cq^2$ and
$\|\ell_{jk}(X)\|_{\mathcal L^q}\le Cq$ for $q\ge2$.  Hence the explicit L\'evy rate following Proposition~\ref{thm:all-rank-chaos} applies with
$(\alpha_h,\alpha_\ell)=(2,1)$, and
$\log p=o(n^{1/13})$ is sufficient for the L\'evy approximation.
\end{proposition}

{This application satisfies Assumption~\ref{ass:block-PSD} with
$\mathcal H_{jk}=\R$ and $\Phi_{jk}(x)=x_jx_k$.  Under the null,
\[
 \mathcal C_{(j,k),(j'k')}=\Gamma_{jk,j'k'},
 \qquad
 \widehat{\mathcal C}_{(j,k),(j'k')}=\widehat\Gamma_{jk,j'k'},
\quad
 \Delta
 =\max_{(j,k),(j'k')\in\mathcal E}
 |\widehat\Gamma_{jk,j'k'}-\Gamma_{jk,j'k'}|.
\]
}

For this application, Step~1 of
Algorithm~\ref{alg:gram-bootstrap} uses the observable scale
\[
 \widehat a_{jk}=\frac1n\sum_{i=1}^n X_{ij}^2X_{ik}^2,
 \qquad
 \widehat W_{n,jk}^{\rm edge}
 =\sqrt{\frac{2}{n(n-1)}}
 \frac{\sum_{i<\ell}X_{ij}X_{ik}X_{\ell j}X_{\ell k}}
 {\widehat a_{jk}}.
\]
Let $\overline{X_jX_k}=n^{-1}\sum_iX_{ij}X_{ik}$ and
\[
 \widehat s_{jk}^2
 =\frac1n\sum_{i=1}^n\left\{X_{ij}X_{ik}-\overline{X_jX_k}\right\}^2.
\]
Define the empirical correlation matrix of the centered edge scores by
\begin{equation*}
 \widehat\Gamma_{jk,j'k'}
 =\frac{n^{-1}\sum_{i=1}^n
       \left\{X_{ij}X_{ik}-\overline{X_jX_k}\right\}
       \left\{X_{ij'}X_{ik'}-\overline{X_{j'}X_{k'}}\right\}}
 {\widehat s_{jk}\widehat s_{j'k'}}.
\end{equation*}
Define
\(
 \bm u_{jk}=\left(X_{1j}X_{1k},\ldots,X_{nj}X_{nk}\right)^\top\),
and \( \mathbf B_{jk}=\mathbf P_n\bm u_{jk}\bm u_{jk}^\top\mathbf P_n.
\)
Then, on the event
$\min_{(j,k)\in\mathcal E}\widehat s_{jk}>0$, the Gram multiplier coordinate
equals $(\widehat Z_{jk}^2-1)/\sqrt2$, where, conditionally on the data,
$\widehat{\bm Z}\sim N(0,\widehat{\bm\Gamma})$.

\begin{theorem}
\label{thm:gaussian-covariance-bootstrap}
Assume $0<c\le\sigma_{jj}\le C<\infty$ uniformly in $j$ and
$H_{0,\mathcal E}$. {For $0<\eta<1/2$, let
\begin{equation*}
 v_n(\eta)=C\left\{
 \sqrt{\frac{\log(8p/\eta)}n}+\frac{\log^3(8p/\eta)}n
 \right\}.
\end{equation*}}
With probability at least $1-\eta$,
\begin{equation}
 \max_{(j,k)\in\mathcal E}
 \left|\frac{\widehat a_{jk}}{\E(X_j^2X_k^2)}-1\right|
 +\max_{(j,k)\in\mathcal E}
 \left|\frac{\widehat s_{jk}^2}{\E(X_j^2X_k^2)}-1\right|
 +\max_{\substack{(j,k)\in\mathcal E\\(j',k')\in\mathcal E}}
 |\widehat\Gamma_{jk,j'k'}-\Gamma_{jk,j'k'}|
 \le v_n(\eta),
 \label{eq:gaussian-score-uniform-concentration}
\end{equation}
where $\bm\Gamma$ and $\widehat{\bm\Gamma}$ are the population and empirical
correlation matrices of the edge scores.  On this event, if $v_n(\eta)\le1/2$,
{Theorem~\ref{thm:gram-bootstrap} applies with
$\Delta\le v_n(\eta)$.

Additionally, the rank-one representation in this application further sharpens the
general Kolmogorov conclusion to}
\begin{equation}
 \sup_t\left|
 \Pp^\#\!\left(\max_{(j,k)\in\mathcal E} Q_{n,jk}^{\#}\le t\right)
 -\Pp\!\left(\max_{(j,k)\in\mathcal E} Q_{jk}\le t\right)
 \right|
 \le C v_n(\eta)^{1/3}
 \{1\vee\log(p/v_n(\eta))\}^{2/3}.
 \label{eq:rank-one-bootstrap-comparison}
\end{equation}
\end{theorem}
}

For any fixed $\alpha$, along any sequence $\eta=\eta_n\downarrow0$,
combining this bound with Theorem~\ref{thm:structured-anticoncentration} and
\eqref{eq:gaussian-score-uniform-concentration}    yields asymptotically valid conditional $(1-\alpha)$ Gram
critical values for
$\max_{(j,k)\in\mathcal E}\widehat W_{n,jk}^{\rm edge}$ whenever the approximation errors vanish.


\subsection{Two-sample distributional testing via  maximum MMD}
\label{sec:mmd-app}
Let $\bm X_1,\ldots,\bm X_n\sim\mathbb P$ and $\bm Y_1,\ldots,\bm Y_n\sim\mathbb Q$ be independent samples
on a high-dimensional product space.  Let
$\mathcal G_1,\ldots,\mathcal G_p$ be prespecified coordinates or scientific
blocks, and {use the same bounded measurable characteristic
kernel $k$ for every block, with separable RKHS and measurable canonical
feature map $\varphi$.}
Define, for
$\bm Z_i=(\bm X_i^\top,\bm Y_i^\top)^\top$,
\begin{equation}
 \Phi_j(\bm Z_i)=\varphi(\bm X_{i,\mathcal G_j})-
             \varphi(\bm Y_{i,\mathcal G_j}),
 \qquad
 h_j(\bm Z_i,\bm Z_k)=\langle\Phi_j(\bm Z_i),\Phi_j(\bm Z_k)\rangle.
 \label{eq:block-MMD-kernel}
\end{equation}
  The simultaneous null
is
$H_0:\mathbb P_{\mathcal G_j}=\mathbb Q_{\mathcal G_j}$ for all $1\le j\le p$.

\begin{proposition}
\label{prop:mmd-reduction}
Under $H_0$, the kernels in \eqref{eq:block-MMD-kernel} are canonical and
positive semidefinite, and
\begin{equation*}
 \E h_j(Z,Z')=
 \|\mu_{\mathbb P,\mathcal G_j}-\mu_{\mathbb Q,\mathcal G_j}\|_{\mathcal H_k}^2
 =\operatorname{MMD}_k^2(\mathbb P_{\mathcal G_j},\mathbb Q_{\mathcal G_j}).
\end{equation*}
{For example, the Gaussian RBF choice
$k(x,y)=\exp\{-\|x-y\|^2/(2\sigma^2)\}$, $\sigma>0$, satisfies $k(x,x)=1$,
so the bounded-diagonal condition below holds with $K_0=1$.}
If $\sup_x k(x,x)\le K_0^2$ and
$\min_j\|h_j\|_{\mathcal L^2(\mathbb P_Z^2)}\ge a_0>0$, then the normalized kernels satisfy
Assumption~\ref{ass:kernel-leverage} with $(\alpha_h,\alpha_\ell)=(0,0)$ and bounded
$\mathcal K$.  Therefore the  L\'evy rate following Proposition~\ref{thm:all-rank-chaos} applies with
$(\alpha_h,\alpha_\ell)=(0,0)$ and
$\log p=o(n^{1/10})$ is sufficient for the L\'evy approximation.
\end{proposition}

The signed target reduces to the
positive-spectrum chaos determined by the covariance operator of
$\Phi_j(Z)$, which is generally infinite rank.
The next corollary  verifies the geometric spectral condition required by
Theorem~\ref{thm:gram-bootstrap}(iii).

\begin{corollary}
\label{cor:gaussian-RBF-Gram}
Consider one-dimensional blocks and the Gaussian RBF kernel
\(
 k_\sigma(x,y)=\exp\{-(x-y)^2/(2\sigma^2)\},\) for some \(\sigma>0.
\)
Under the  null, suppose
\[
 \mathbb P_{\mathcal G_j}=\mathbb Q_{\mathcal G_j}=N(m_j,s_j^2),
 \qquad 0<s_-\le s_j\le s_+<\infty,
\]
uniformly in $j$.  Then Assumption~\ref{ass:block-PSD} holds, and there exist constants
$c_1,c_2>0$ and $\beta\ge\alpha>0$, depending only on
$s_-,s_+$ and $\sigma$, such that the eigenvalues in
\eqref{eq:PSD-Gaussian-block-representation} satisfy
\[
 c_2e^{-\beta r}\le\mu_{jr}\le c_1e^{-\alpha r},
 \qquad 1\le j\le p,\ r\ge1.
\]
 For
$0<\eta<1/2$, with probability at least $1-\eta$,
\[
 \Delta\lesssim
 \sqrt{\frac{\log(8p^2/\eta)}n}
 +\frac{\log(8p^2/\eta)}n.
\]
On this event, the conclusion of Theorem \ref{thm:gram-bootstrap}(iii) holds.
\end{corollary}

The quadratic-time MMD representation of \citet{GrettonEtAl2012} therefore extends naturally to simultaneous testing across multiple blocks. This setting differs from normal approximations for a single high-dimensional studentized MMD \cite{GaoShao2023}: here, \(p\) indexes distinct blockwise hypotheses, and coordinates with low effective rank need not Gaussianize.



\section{Numerical simulations}\label{sec:simulations}

\subsection{Simulation setup}

We first examine how the spectrum affects the accuracy of Gaussian and
chaos calibration, then study size and power in the two applications of
Section~\ref{sec:applications}.
Throughout, ``Oracle'' denotes joint population-chaos calibration and
``Prop.'' the proposed Gaussian multiplier bootstrap.
The application experiments vary dependence
while keeping the population marginal distributions fixed, so that the comparison with the
oracle and Bonferroni baseline measures the benefit of joint calibration.
Rejection frequencies are based on 500 Monte Carlo replications.
We report level-0.05 results here. The supplementary material provides the
corresponding level-0.10 results, implementation details, and additional
experimental designs.


\subsection{Spectral phase transition}
\label{sec:simulation-mechanism}

Let \(\bm Y_i\in\R^{64}\) be independent with mean zero and identity covariance.
For \(j\le p\), define
\[
 h_j(\bm Y_i,\bm Y_k)=\bm Y_i^\top\mathbf A_j\bm Y_k,\qquad
 \mathbf A_j=\mathbf R_j\mathbf D_r\mathbf R_j^\top,\qquad \|\mathbf D_r\|_{\mathrm F}=1,
\]
where \(\mathbf R_j\) is orthogonal.  Four spectra are used:
\[
 \mathbf D_1=(1),\qquad
 \mathbf D_3\propto\operatorname{diag}(1,-0.8,0.6),\qquad
 \mathbf D_r=r^{-1/2}\operatorname{diag}(1,-1,\ldots),\quad r\in\{16,64\}.
\]
Zeros pad each diagonal to dimension \(64\).  The flat signed spectra have
\(r_4=r\), while the signed rank-three spectrum has
\(r_4=4/1.5392\approx2.6\).  The four choices therefore range from a single dominant eigenvalue to
a flat spectrum, including both positive and signed kernels.

We vary dependence without changing these spectra by orthogonalizing
\[
 \sqrt{\rho_{\rm rot}}\mathbf I_{64}
 +\sqrt{1-\rho_{\rm rot}}\mathbf G_j/\sqrt{64},
\]
where \(\mathbf G_j\) has independent standard Gaussian entries and
\(\rho_{\rm rot}\in\{0,0.95\}\). The cross-kernel covariance is
\(\tr(\mathbf A_j\mathbf A_k)\); thus \(\rho_{\rm rot}\) controls the
rotations rather than specifying their correlations. Gaussian scores are
used with \(n\in\{200,500\}\) and \(p\in\{200,1000\}\). We also consider
Rademacher scores at \(n=500\) and \(p=1000\) under both rotation choices.
The proposed bootstrap is compared with the oracle joint signed-chaos law
and the covariance-matched Gaussian maximum.

Table~\ref{tab:mechanism-primary} reports rejection frequencies at level 0.05.
The proposed rejection frequencies remain close to those of the oracle,
whereas ordinary Gaussian calibration substantially overrejects for the
rank-one and signed rank-three kernels. This overrejection decreases as the
effective rank increases to 64, consistent with the Gaussianization result
in Theorem~\ref{thm:chaos-Gaussianization}.

\begin{table}[!htbp]
\centering
\caption{Rejection frequencies at level 0.05 for each spectral setting.  Gaussian denotes the covariance-matched Gaussian approximation.}
\label{tab:mechanism-primary}
\begingroup\footnotesize
\setlength{\tabcolsep}{0.8pt}
\renewcommand{\arraystretch}{1.03}
\begin{tabular*}{\textwidth}{@{\extracolsep{\fill}}llllrrrrrr@{}}
\toprule
 &  &  &  & \multicolumn{3}{c}{$\rho_{\rm rot}=0$} & \multicolumn{3}{c}{$\rho_{\rm rot}=0.95$} \\
\cmidrule(lr){5-7}\cmidrule(lr){8-10}
$n$ & $p$ & Scores & $r_4$ & Oracle & Prop. & Gaussian & Oracle & Prop. & Gaussian \\
\midrule
200 & 200 & Gaussian & 1 & 0.046 & 0.042 & 0.886 & 0.078 & 0.062 & 0.786 \\
200 & 200 & Gaussian & 2.6 & 0.038 & 0.044 & 0.622 & 0.074 & 0.056 & 0.494 \\
200 & 200 & Gaussian & 16 & 0.054 & 0.038 & 0.154 & 0.066 & 0.042 & 0.174 \\
200 & 200 & Gaussian & 64 & 0.060 & 0.056 & 0.074 & 0.034 & 0.038 & 0.068 \\
\addlinespace[3pt]
200 & 1000 & Gaussian & 1 & 0.046 & 0.046 & 0.988 & 0.060 & 0.042 & 0.870 \\
200 & 1000 & Gaussian & 2.6 & 0.062 & 0.040 & 0.870 & 0.050 & 0.038 & 0.588 \\
200 & 1000 & Gaussian & 16 & 0.038 & 0.046 & 0.284 & 0.044 & 0.050 & 0.198 \\
200 & 1000 & Gaussian & 64 & 0.052 & 0.048 & 0.084 & 0.038 & 0.034 & 0.076 \\
\addlinespace[3pt]
500 & 200 & Gaussian & 1 & 0.046 & 0.046 & 0.866 & 0.076 & 0.056 & 0.754 \\
500 & 200 & Gaussian & 2.6 & 0.054 & 0.052 & 0.624 & 0.048 & 0.050 & 0.520 \\
500 & 200 & Gaussian & 16 & 0.036 & 0.034 & 0.160 & 0.042 & 0.046 & 0.116 \\
500 & 200 & Gaussian & 64 & 0.052 & 0.050 & 0.072 & 0.052 & 0.046 & 0.054 \\
\addlinespace[3pt]
500 & 1000 & Gaussian & 1 & 0.042 & 0.038 & 0.982 & 0.044 & 0.038 & 0.890 \\
500 & 1000 & Gaussian & 2.6 & 0.036 & 0.040 & 0.854 & 0.038 & 0.034 & 0.608 \\
500 & 1000 & Gaussian & 16 & 0.040 & 0.036 & 0.274 & 0.052 & 0.046 & 0.172 \\
500 & 1000 & Gaussian & 64 & 0.052 & 0.034 & 0.046 & 0.058 & 0.044 & 0.084 \\
\addlinespace[3pt]
500 & 1000 & Rademacher & 1 & 0.044 & 0.042 & 0.974 & 0.026 & 0.046 & 0.830 \\
500 & 1000 & Rademacher & 2.6 & 0.042 & 0.046 & 0.878 & 0.058 & 0.048 & 0.580 \\
500 & 1000 & Rademacher & 16 & 0.042 & 0.038 & 0.224 & 0.070 & 0.046 & 0.194 \\
500 & 1000 & Rademacher & 64 & 0.062 & 0.048 & 0.070 & 0.038 & 0.034 & 0.064 \\
\bottomrule
\end{tabular*}
\endgroup
\end{table}

\begin{figure}[!htbp]
\centering
\includegraphics[width=0.85\textwidth]{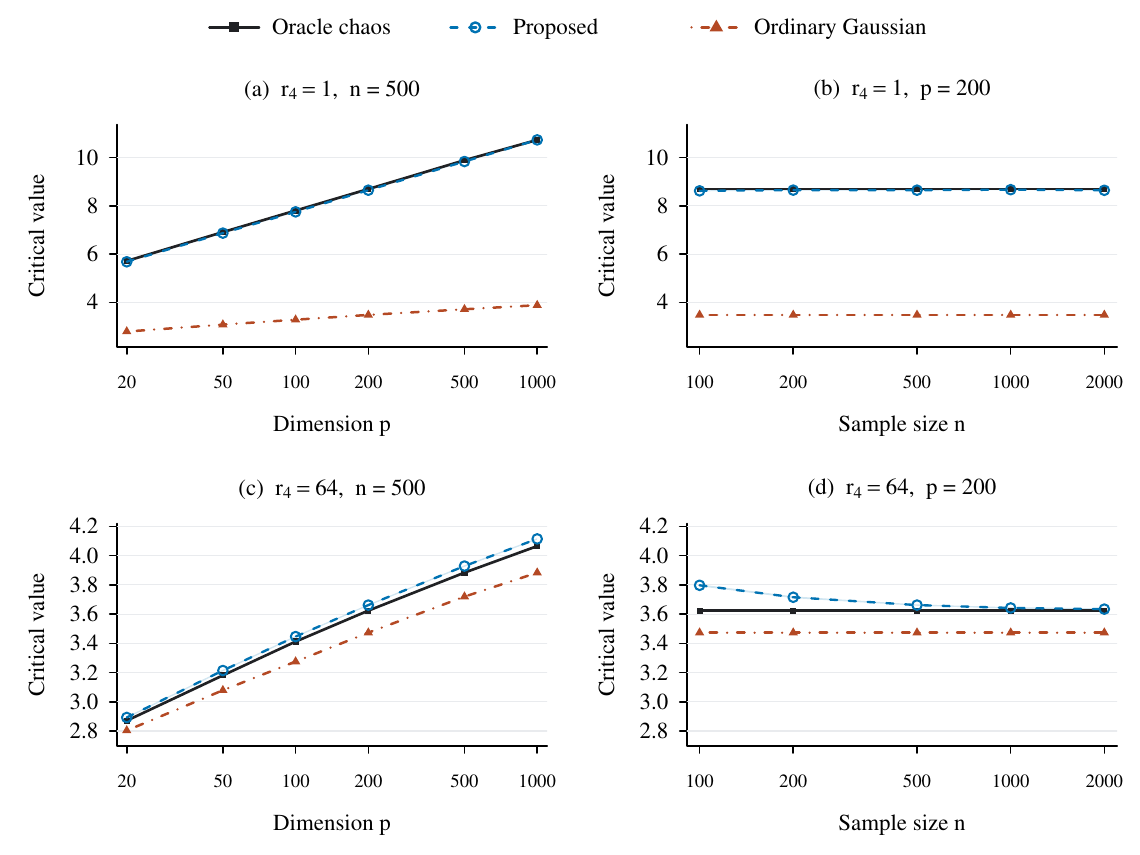}
\caption{Oracle, proposed, and ordinary Gaussian critical values.
Shading gives pointwise 95\% Monte Carlo intervals for the mean proposed
critical value.}
\label{fig:critical-value-transition}
\end{figure}

Figure~\ref{fig:critical-value-transition} shows the effects of dimension
and sample size on the 95\% critical value, using Gaussian scores and
$\rho_{\rm rot}=0$. At $r_4=1$, the proposed curve closely follows the
population chaos quantile, while the Gaussian approximation is too small.
At $r_4=64$, the two population quantiles are much closer; the proposed
critical value is slightly high for small $n$ and approaches the oracle
as $n$ grows. Increasing $n$ reduces bootstrap estimation error but does
not close the Gaussian-to-chaos gap for a fixed spectrum and dimension.

\subsection{Gaussian covariance-graph testing}
We study how dependence among tested edges affects the covariance test in
Section~\ref{sec:gaussian-covariance-app}. Let $U,F,e_1,\ldots,e_p$ be
independent standard Gaussian variables and set
$V_j=\sqrt\rho F+\sqrt{1-\rho}e_j$. The $p$ candidate edges join $U$ to
the $V_j$'s. Their covariances are zero under the null, but the edge scores
have correlation $\rho$ and the corresponding chaos coordinates have
correlation $\rho^2$. We use $n\in\{200,500\}$ paired differences,
$p\in\{50,200\}$, and $\rho\in\{0,0.2,0.5,0.8\}$. Under the alternatives,
the first $s\in\{1,5\}$ variables are replaced by $V_j+\delta U$, where
$\delta=c\sqrt{\log(p)/n}$ and $c\in\{0.75,1.25,1.75\}$.

The oracle uses the joint population chaos, while Bonferroni uses its
centered chi-square marginal at tail probability $\alpha/p$. Both are
applied to the same statistic as the proposed bootstrap. At $\rho=0$,
the population chaos coordinates are independent, although the
finite-sample scores share $U$.

The proposed test closely tracks the oracle in both size and power; see
 Tables~\ref{tab:cov-null-primary}--\ref{tab:cov-power-primary}. At
$\rho=0$, Bonferroni is also close to the joint procedures. As dependence
increases, it becomes conservative and loses power, whereas the bootstrap
continues to reproduce the oracle calibration. For example, at
$n=500$, $p=200$, $s=5$, $\rho=0.8$, and $c=1.25$, proposed and oracle
power are 0.622 and 0.632, compared with 0.336 for Bonferroni.
This gain comes from using the joint
distribution of the edge scores: Bonferroni assigns $\alpha/p$ to each
test.

We also consider diagonal covariance matrices, banded matrices with
first off-diagonal 0.30, equicorrelated blocks of size 10 with correlation
0.70, and a sparse cycle-and-matching structure with correlation 0.20.
We take $d\in\{50,100\}$ and test either the first $d$ null edges in
lexicographic order or all null edges. The proposed and oracle procedures
remain close, although both can be conservative when many edges are tested; see
Table~\ref{tab:covariance-structures}. Ordinary Gaussian calibration
substantially overrejects, as expected from the rank-one representation in
Proposition~\ref{prop:gaussian-covariance-reduction}. This contrast shows
why retaining the quadratic-chaos law matters even with Gaussian data.
For diagonal covariance matrices with all edges tested, we also include
the test of \citet{CaiJiang2011}, denoted by C--J, whose rejection
frequencies are similarly conservative. The supplementary material
reports the level-0.10 results for these same designs and considers a
separate sparse-edge alternative with all
$\binom d2$ edges tested, increasing the number of hypotheses to 1225 or 4950.

\begin{table}[!htbp]
\centering
\caption{Null rejection frequencies for covariance testing at level 0.05. Bonf. uses population marginal-chaos quantiles.}
\label{tab:cov-null-primary}
\begingroup\footnotesize
\setlength{\tabcolsep}{0.8pt}
\renewcommand{\arraystretch}{1.08}
\begin{tabular*}{\textwidth}{@{\extracolsep{\fill}}llrrrrrrrrrrrr@{}}
\toprule
 &  & \multicolumn{3}{c}{$\rho=0$} & \multicolumn{3}{c}{$\rho=0.2$} & \multicolumn{3}{c}{$\rho=0.5$} & \multicolumn{3}{c}{$\rho=0.8$} \\
\cmidrule(lr){3-5}\cmidrule(lr){6-8}\cmidrule(lr){9-11}\cmidrule(lr){12-14}
$n$ & $p$ & Oracle & Prop. & Bonf. & Oracle & Prop. & Bonf. & Oracle & Prop. & Bonf. & Oracle & Prop. & Bonf. \\
\midrule
200 & 50 & 0.048 & 0.046 & 0.044 & 0.042 & 0.060 & 0.038 & 0.050 & 0.050 & 0.028 & 0.050 & 0.044 & 0.014 \\
\addlinespace[3pt]
200 & 200 & 0.026 & 0.028 & 0.024 & 0.024 & 0.026 & 0.020 & 0.034 & 0.036 & 0.014 & 0.036 & 0.032 & 0.006 \\
\addlinespace[3pt]
500 & 50 & 0.050 & 0.046 & 0.046 & 0.050 & 0.050 & 0.044 & 0.032 & 0.036 & 0.018 & 0.064 & 0.064 & 0.010 \\
\addlinespace[3pt]
500 & 200 & 0.058 & 0.046 & 0.056 & 0.046 & 0.042 & 0.036 & 0.036 & 0.030 & 0.018 & 0.050 & 0.050 & 0.008 \\
\bottomrule
\end{tabular*}
\endgroup
\par\vspace{6pt}
\caption{Global power for covariance testing at level 0.05. Bonf. uses population marginal-chaos quantiles.}
\label{tab:cov-power-primary}
\begingroup\footnotesize
\setlength{\tabcolsep}{0.8pt}
\renewcommand{\arraystretch}{1.00}
\begin{tabular*}{\textwidth}{@{\extracolsep{\fill}}llllrrrrrrrrrrrr@{}}
\toprule
 &  &  &  & \multicolumn{3}{c}{$\rho=0$} & \multicolumn{3}{c}{$\rho=0.2$} & \multicolumn{3}{c}{$\rho=0.5$} & \multicolumn{3}{c}{$\rho=0.8$} \\
\cmidrule(lr){5-7}\cmidrule(lr){8-10}\cmidrule(lr){11-13}\cmidrule(lr){14-16}
$n$ & $p$ & $s$ & $c$ & Oracle & Prop. & Bonf. & Oracle & Prop. & Bonf. & Oracle & Prop. & Bonf. & Oracle & Prop. & Bonf. \\
\midrule
200 & 50 & 1 & 0.75 & 0.054 & 0.056 & 0.048 & 0.060 & 0.056 & 0.050 & 0.076 & 0.078 & 0.052 & 0.104 & 0.112 & 0.026 \\
200 & 50 & 1 & 1.25 & 0.178 & 0.176 & 0.176 & 0.192 & 0.200 & 0.188 & 0.226 & 0.214 & 0.152 & 0.356 & 0.358 & 0.158 \\
200 & 50 & 1 & 1.75 & 0.468 & 0.472 & 0.466 & 0.488 & 0.496 & 0.476 & 0.564 & 0.564 & 0.490 & 0.716 & 0.728 & 0.490 \\
200 & 50 & 5 & 0.75 & 0.154 & 0.150 & 0.148 & 0.128 & 0.134 & 0.122 & 0.148 & 0.144 & 0.104 & 0.180 & 0.178 & 0.058 \\
200 & 50 & 5 & 1.25 & 0.592 & 0.602 & 0.590 & 0.542 & 0.556 & 0.520 & 0.530 & 0.528 & 0.454 & 0.552 & 0.542 & 0.298 \\
200 & 50 & 5 & 1.75 & 0.958 & 0.958 & 0.958 & 0.904 & 0.906 & 0.900 & 0.870 & 0.870 & 0.800 & 0.876 & 0.874 & 0.664 \\
\addlinespace[3pt]
200 & 200 & 1 & 0.75 & 0.048 & 0.052 & 0.046 & 0.044 & 0.042 & 0.036 & 0.060 & 0.056 & 0.026 & 0.106 & 0.104 & 0.022 \\
200 & 200 & 1 & 1.25 & 0.174 & 0.174 & 0.166 & 0.152 & 0.168 & 0.144 & 0.220 & 0.234 & 0.150 & 0.428 & 0.434 & 0.168 \\
200 & 200 & 1 & 1.75 & 0.524 & 0.526 & 0.518 & 0.538 & 0.536 & 0.518 & 0.598 & 0.608 & 0.492 & 0.786 & 0.798 & 0.498 \\
200 & 200 & 5 & 0.75 & 0.112 & 0.124 & 0.112 & 0.098 & 0.098 & 0.082 & 0.104 & 0.102 & 0.062 & 0.184 & 0.194 & 0.054 \\
200 & 200 & 5 & 1.25 & 0.568 & 0.580 & 0.562 & 0.508 & 0.522 & 0.494 & 0.534 & 0.548 & 0.406 & 0.640 & 0.642 & 0.280 \\
200 & 200 & 5 & 1.75 & 0.958 & 0.960 & 0.956 & 0.962 & 0.964 & 0.946 & 0.892 & 0.896 & 0.814 & 0.926 & 0.928 & 0.724 \\
\addlinespace[3pt]
500 & 50 & 1 & 0.75 & 0.066 & 0.066 & 0.062 & 0.080 & 0.080 & 0.074 & 0.072 & 0.070 & 0.050 & 0.136 & 0.136 & 0.060 \\
500 & 50 & 1 & 1.25 & 0.236 & 0.234 & 0.232 & 0.234 & 0.242 & 0.224 & 0.272 & 0.272 & 0.208 & 0.404 & 0.398 & 0.212 \\
500 & 50 & 1 & 1.75 & 0.564 & 0.564 & 0.562 & 0.558 & 0.566 & 0.550 & 0.588 & 0.584 & 0.510 & 0.764 & 0.774 & 0.550 \\
500 & 50 & 5 & 0.75 & 0.168 & 0.166 & 0.166 & 0.204 & 0.194 & 0.190 & 0.166 & 0.164 & 0.108 & 0.220 & 0.210 & 0.074 \\
500 & 50 & 5 & 1.25 & 0.672 & 0.672 & 0.664 & 0.592 & 0.602 & 0.580 & 0.594 & 0.582 & 0.480 & 0.566 & 0.562 & 0.332 \\
500 & 50 & 5 & 1.75 & 0.978 & 0.980 & 0.976 & 0.958 & 0.952 & 0.950 & 0.908 & 0.906 & 0.862 & 0.912 & 0.906 & 0.752 \\
\addlinespace[3pt]
500 & 200 & 1 & 0.75 & 0.072 & 0.070 & 0.070 & 0.068 & 0.072 & 0.062 & 0.094 & 0.088 & 0.040 & 0.142 & 0.148 & 0.028 \\
500 & 200 & 1 & 1.25 & 0.228 & 0.224 & 0.224 & 0.220 & 0.216 & 0.196 & 0.288 & 0.280 & 0.206 & 0.416 & 0.416 & 0.166 \\
500 & 200 & 1 & 1.75 & 0.606 & 0.610 & 0.602 & 0.638 & 0.634 & 0.620 & 0.682 & 0.678 & 0.574 & 0.830 & 0.832 & 0.588 \\
500 & 200 & 5 & 0.75 & 0.150 & 0.138 & 0.146 & 0.140 & 0.132 & 0.118 & 0.160 & 0.158 & 0.102 & 0.198 & 0.188 & 0.044 \\
500 & 200 & 5 & 1.25 & 0.652 & 0.658 & 0.646 & 0.602 & 0.590 & 0.564 & 0.600 & 0.594 & 0.474 & 0.632 & 0.622 & 0.336 \\
500 & 200 & 5 & 1.75 & 0.994 & 0.992 & 0.992 & 0.974 & 0.980 & 0.964 & 0.956 & 0.958 & 0.896 & 0.958 & 0.956 & 0.790 \\
\bottomrule
\end{tabular*}
\endgroup
\end{table}

\begin{table}[ht]
\centering
\caption{Null rejection frequencies at level 0.05 under different covariance structures. The edge sets comprise either all null edges or the first $d$ null edges in lexicographic order. C--J denotes the coherence test of \citet{CaiJiang2011}; dashes indicate settings outside its comparison regime. Bonf. uses population marginal-chaos quantiles.}
\label{tab:covariance-structures}
\begingroup\footnotesize
\setlength{\tabcolsep}{0.8pt}
\renewcommand{\arraystretch}{1.08}
\begin{tabular*}{\textwidth}{@{\extracolsep{\fill}}lllrrrrrrrrrr@{}}
\toprule
 &  &  & \multicolumn{5}{c}{$n=200$} & \multicolumn{5}{c}{$n=500$} \\
\cmidrule(lr){4-8}\cmidrule(lr){9-13}
$d$ & Covariance & Edges & Oracle & Prop. & Bonf. & C--J & Gauss. & Oracle & Prop. & Bonf. & C--J & Gauss. \\
\midrule
50 & Banded & All & 0.024 & 0.028 & 0.020 & -- & 1.000 & 0.010 & 0.016 & 0.016 & -- & 1.000 \\
50 & Banded & First $d$ & 0.040 & 0.036 & 0.032 & -- & 0.588 & 0.060 & 0.056 & 0.052 & -- & 0.596 \\
\addlinespace[3pt]
50 & Equicorrelated blocks & All & 0.042 & 0.032 & 0.018 & -- & 0.946 & 0.056 & 0.032 & 0.014 & -- & 0.940 \\
50 & Equicorrelated blocks & First $d$ & 0.040 & 0.046 & 0.022 & -- & 0.374 & 0.046 & 0.050 & 0.032 & -- & 0.376 \\
\addlinespace[3pt]
50 & Diagonal & All & 0.018 & 0.018 & 0.018 & 0.014 & 1.000 & 0.030 & 0.038 & 0.032 & 0.030 & 1.000 \\
50 & Diagonal & First $d$ & 0.040 & 0.052 & 0.046 & -- & 0.658 & 0.044 & 0.038 & 0.038 & -- & 0.636 \\
\addlinespace[3pt]
50 & Sparse cycle/matching & All & 0.018 & 0.022 & 0.018 & -- & 1.000 & 0.038 & 0.040 & 0.036 & -- & 1.000 \\
50 & Sparse cycle/matching & First $d$ & 0.022 & 0.036 & 0.030 & -- & 0.640 & 0.046 & 0.040 & 0.042 & -- & 0.606 \\
\addlinespace[3pt]
100 & Banded & First $d$ & 0.024 & 0.024 & 0.024 & -- & 0.830 & 0.050 & 0.048 & 0.042 & -- & 0.820 \\
\addlinespace[3pt]
100 & Equicorrelated blocks & First $d$ & 0.036 & 0.046 & 0.024 & -- & 0.622 & 0.054 & 0.054 & 0.036 & -- & 0.622 \\
\addlinespace[3pt]
100 & Diagonal & All & 0.014 & 0.018 & 0.016 & 0.014 & 1.000 & 0.030 & 0.026 & 0.028 & 0.024 & 1.000 \\
100 & Diagonal & First $d$ & 0.024 & 0.032 & 0.028 & -- & 0.822 & 0.052 & 0.044 & 0.044 & -- & 0.824 \\
\addlinespace[3pt]
100 & Sparse cycle/matching & First $d$ & 0.026 & 0.026 & 0.028 & -- & 0.782 & 0.030 & 0.044 & 0.040 & -- & 0.820 \\
\bottomrule
\end{tabular*}
\endgroup
\end{table}

\FloatBarrier

\subsection{Two-sample testing via maximum MMD}
We next consider the blockwise two-sample test in Section~\ref{sec:mmd-app}.
Within each sample, let $Z_j=\sqrt\rho F+\sqrt{1-\rho}e_j$ with independent
standard Gaussian $F,e_1,\ldots,e_p$, and transform the coordinates to
centered, variance-standardized $t_3$ or log-normal marginals. Observations
and the two samples are independent. We use $n\in\{100,200\}$,
$p=20$ univariate blocks, and $\rho\in\{0,0.2,0.5,0.8\}$; here $\rho$
is the latent Gaussian correlation. Each block uses an RBF kernel with
bandwidth 1. Alternatives shift the first $s\in\{1,5\}$ coordinates of
the second sample by 0.15, 0.30, or 0.45.

As in the covariance experiment, the oracle and Bonferroni use the same
population marginal chaos distributions.  The proposed method  estimates the joint distribution
from the centered Gram matrices using Gaussian multipliers. The population-reference construction and calibration details are
given in the supplementary material.

The proposed bootstrap and oracle give similar rejection frequencies across
the dependence settings; see Table~\ref{tab:mmd-null-primary}. Calibration
is less accurate in small log-normal samples: at $n=100$ and $\rho=0.5$,
both reject with frequency 0.094 at nominal level 0.05. For $n=200$,
proposed rejection frequencies range from 0.044 to 0.064 across the
individual settings. Replacing the estimated scale by its population value on the same
500 datasets reduces oracle rejection from 0.094 to 0.068 in the
log-normal setting with $n=100$ and $\rho=0.5$. Thus scale estimation
contributes to the small-sample deviation; increasing $n$ improves
calibration for both joint procedures.

Dependence has the clearest effect on power at the smaller shifts; see Table~\ref{tab:mmd-power-primary}. With $t_3$ marginals, $n=200$, $s=1$,
a shift of 0.15, and $\rho=0.8$, proposed and oracle power are 0.310 and
0.308, compared with 0.188 for Bonferroni. The difference is small at
$\rho=0$, and all three procedures
approach unit power at the largest shift. Thus joint calibration improves
sensitivity where a marginal correction discards useful dependence
information, without requiring the population eigensystem.
The supplementary material gives the corresponding level-0.10 results
and considers additional independent-block designs, with
other marginal laws, block dimensions, bandwidths, and distributional
changes, and compare the proposed procedure with permutation.


\begin{table}[ht]
\centering
\caption{Null rejection frequencies for maximum MMD at level 0.05. Each setting has $p=20$ blocks and bandwidth 1. Bonf. uses population marginal-chaos quantiles.}
\label{tab:mmd-null-primary}
\begingroup\footnotesize
\setlength{\tabcolsep}{0.8pt}
\renewcommand{\arraystretch}{1.08}
\begin{tabular*}{\textwidth}{@{\extracolsep{\fill}}llrrrrrrrrrrrr@{}}
\toprule
 &  & \multicolumn{3}{c}{$\rho=0$} & \multicolumn{3}{c}{$\rho=0.2$} & \multicolumn{3}{c}{$\rho=0.5$} & \multicolumn{3}{c}{$\rho=0.8$} \\
\cmidrule(lr){3-5}\cmidrule(lr){6-8}\cmidrule(lr){9-11}\cmidrule(lr){12-14}
$n$ & Marginal law & Oracle & Prop. & Bonf. & Oracle & Prop. & Bonf. & Oracle & Prop. & Bonf. & Oracle & Prop. & Bonf. \\
\midrule
100 & $t_3$ & 0.050 & 0.056 & 0.050 & 0.050 & 0.050 & 0.048 & 0.044 & 0.050 & 0.034 & 0.062 & 0.056 & 0.032 \\
\addlinespace[3pt]
100 & Log-normal & 0.058 & 0.066 & 0.056 & 0.062 & 0.064 & 0.060 & 0.094 & 0.094 & 0.086 & 0.060 & 0.062 & 0.036 \\
\addlinespace[3pt]
200 & $t_3$ & 0.066 & 0.064 & 0.064 & 0.050 & 0.058 & 0.050 & 0.058 & 0.056 & 0.044 & 0.040 & 0.044 & 0.028 \\
\addlinespace[3pt]
200 & Log-normal & 0.052 & 0.048 & 0.050 & 0.044 & 0.044 & 0.042 & 0.066 & 0.064 & 0.056 & 0.052 & 0.050 & 0.030 \\
\bottomrule
\end{tabular*}
\endgroup
\end{table}

\begin{table}[ht]
\centering
\caption{Global power for maximum MMD at level 0.05. Each setting has $p=20$ blocks and bandwidth 1. Bonf. uses population marginal-chaos quantiles.}
\label{tab:mmd-power-primary}
\begingroup\footnotesize
\setlength{\tabcolsep}{0.8pt}
\renewcommand{\arraystretch}{1.08}
\begin{tabular*}{\textwidth}{@{\extracolsep{\fill}}llllrrrrrrrrrrrr@{}}
\toprule
 &  &  &  & \multicolumn{3}{c}{$\rho=0$} & \multicolumn{3}{c}{$\rho=0.2$} & \multicolumn{3}{c}{$\rho=0.5$} & \multicolumn{3}{c}{$\rho=0.8$} \\
\cmidrule(lr){5-7}\cmidrule(lr){8-10}\cmidrule(lr){11-13}\cmidrule(lr){14-16}
$n$ & Marginal law & $s$ & Shift & Oracle & Prop. & Bonf. & Oracle & Prop. & Bonf. & Oracle & Prop. & Bonf. & Oracle & Prop. & Bonf. \\
\midrule
100 & $t_3$ & 1 & 0.15 & 0.116 & 0.116 & 0.116 & 0.126 & 0.126 & 0.126 & 0.128 & 0.132 & 0.108 & 0.116 & 0.124 & 0.076 \\
100 & $t_3$ & 1 & 0.30 & 0.462 & 0.462 & 0.462 & 0.440 & 0.452 & 0.436 & 0.428 & 0.434 & 0.400 & 0.554 & 0.556 & 0.430 \\
100 & $t_3$ & 1 & 0.45 & 0.882 & 0.888 & 0.882 & 0.906 & 0.906 & 0.896 & 0.916 & 0.914 & 0.902 & 0.950 & 0.950 & 0.902 \\
100 & $t_3$ & 5 & 0.15 & 0.302 & 0.302 & 0.296 & 0.254 & 0.258 & 0.246 & 0.252 & 0.248 & 0.218 & 0.232 & 0.238 & 0.152 \\
100 & $t_3$ & 5 & 0.30 & 0.932 & 0.932 & 0.932 & 0.910 & 0.906 & 0.900 & 0.840 & 0.838 & 0.806 & 0.814 & 0.820 & 0.726 \\
100 & $t_3$ & 5 & 0.45 & 1.000 & 1.000 & 1.000 & 1.000 & 1.000 & 1.000 & 1.000 & 1.000 & 0.998 & 0.990 & 0.990 & 0.980 \\
\addlinespace[3pt]
100 & Log-normal & 1 & 0.15 & 0.182 & 0.188 & 0.180 & 0.184 & 0.194 & 0.184 & 0.212 & 0.206 & 0.184 & 0.244 & 0.252 & 0.160 \\
100 & Log-normal & 1 & 0.30 & 0.884 & 0.884 & 0.884 & 0.870 & 0.870 & 0.866 & 0.890 & 0.890 & 0.872 & 0.904 & 0.910 & 0.858 \\
100 & Log-normal & 1 & 0.45 & 1.000 & 1.000 & 1.000 & 1.000 & 1.000 & 1.000 & 1.000 & 1.000 & 1.000 & 1.000 & 1.000 & 1.000 \\
100 & Log-normal & 5 & 0.15 & 0.504 & 0.512 & 0.500 & 0.476 & 0.474 & 0.466 & 0.436 & 0.438 & 0.388 & 0.448 & 0.444 & 0.328 \\
100 & Log-normal & 5 & 0.30 & 1.000 & 1.000 & 1.000 & 0.998 & 0.998 & 0.998 & 0.990 & 0.990 & 0.984 & 0.998 & 1.000 & 0.988 \\
100 & Log-normal & 5 & 0.45 & 1.000 & 1.000 & 1.000 & 1.000 & 1.000 & 1.000 & 1.000 & 1.000 & 1.000 & 1.000 & 1.000 & 1.000 \\
\addlinespace[3pt]
200 & $t_3$ & 1 & 0.15 & 0.194 & 0.206 & 0.194 & 0.188 & 0.198 & 0.184 & 0.210 & 0.202 & 0.184 & 0.308 & 0.310 & 0.188 \\
200 & $t_3$ & 1 & 0.30 & 0.882 & 0.880 & 0.882 & 0.880 & 0.874 & 0.878 & 0.900 & 0.898 & 0.884 & 0.916 & 0.918 & 0.866 \\
200 & $t_3$ & 1 & 0.45 & 0.998 & 0.998 & 0.998 & 1.000 & 1.000 & 1.000 & 1.000 & 1.000 & 1.000 & 1.000 & 1.000 & 1.000 \\
200 & $t_3$ & 5 & 0.15 & 0.636 & 0.642 & 0.636 & 0.546 & 0.548 & 0.542 & 0.506 & 0.514 & 0.456 & 0.448 & 0.448 & 0.342 \\
200 & $t_3$ & 5 & 0.30 & 1.000 & 1.000 & 1.000 & 1.000 & 1.000 & 1.000 & 0.998 & 1.000 & 0.996 & 0.996 & 0.996 & 0.968 \\
200 & $t_3$ & 5 & 0.45 & 1.000 & 1.000 & 1.000 & 1.000 & 1.000 & 1.000 & 1.000 & 1.000 & 1.000 & 1.000 & 1.000 & 1.000 \\
\addlinespace[3pt]
200 & Log-normal & 1 & 0.15 & 0.458 & 0.466 & 0.458 & 0.462 & 0.452 & 0.454 & 0.480 & 0.470 & 0.436 & 0.538 & 0.550 & 0.416 \\
200 & Log-normal & 1 & 0.30 & 0.996 & 0.996 & 0.996 & 0.998 & 0.998 & 0.998 & 0.998 & 0.998 & 0.998 & 1.000 & 1.000 & 1.000 \\
200 & Log-normal & 1 & 0.45 & 1.000 & 1.000 & 1.000 & 1.000 & 1.000 & 1.000 & 1.000 & 1.000 & 1.000 & 1.000 & 1.000 & 1.000 \\
200 & Log-normal & 5 & 0.15 & 0.922 & 0.924 & 0.918 & 0.888 & 0.886 & 0.876 & 0.856 & 0.858 & 0.832 & 0.802 & 0.800 & 0.688 \\
200 & Log-normal & 5 & 0.30 & 1.000 & 1.000 & 1.000 & 1.000 & 1.000 & 1.000 & 1.000 & 1.000 & 1.000 & 1.000 & 1.000 & 1.000 \\
200 & Log-normal & 5 & 0.45 & 1.000 & 1.000 & 1.000 & 1.000 & 1.000 & 1.000 & 1.000 & 1.000 & 1.000 & 1.000 & 1.000 & 1.000 \\
\bottomrule
\end{tabular*}
\endgroup
\end{table}

\FloatBarrier

\section{Concluding remarks}
We have developed a joint Gaussian-chaos approximation for maxima of
canonical order-two $U$-statistics that retains both the quadratic structure
and dependence across kernels. The fourth-order effective rank describes
the transition to Gaussian approximation, while positive spectral structure
yields sharper anti-concentration bounds than those available for general
signed spectra. For positive-semidefinite kernels, the centered-Gram
bootstrap provides joint calibration without estimating population
eigensystems. In the dependent covariance and two-sample experiments, it
closely follows the joint oracle and gains power over Bonferroni calibration.
These results highlight the benefit of accounting for both marginal
non-Gaussianity and dependence in simultaneous inference.

Natural extensions include estimated kernels, higher-order statistics,
and dependent observations. For estimated kernels, sample splitting offers
a conditional approximation route; the main task is to control the induced
first Hoeffding projection and perturbation of the quadratic term. The
present rates provide a benchmark when these errors are of smaller order
after anti-concentration transfer. Higher-order canonical $U$-statistics
call for approximation by higher Wiener chaos and control of contractions
of several orders \cite{NourdinPeccatiReinert2010,DoblerPeccati2019}.
The degree-dependent anti-concentration exponent in
\cite{CarberyWright2001} also changes the loss in passing to Kolmogorov
distance. For short-range dependent observations, dependent wild multipliers
\cite{LeuchtNeumann2013,Chwialkowski2014} offer a way to reproduce long-run
feature covariances. Establishing a joint chaos approximation with this
covariance structure would allow the comparison arguments to be reused,
with rates governed by serial dependence and lag truncation or blocking.



\section*{Supplementary material}
Additional simulation results, auxiliary results, and complete proofs
are provided in the supplementary material.

\section*{Acknowledgment}
\par Qirui Hu's research was supported by National Natural Science Foundation of China (NSFC) (Grant No.~12601520), the Shanghai Engineering Research Center of Finance Intelligence (Grant No.~19DZ2254600) and by TRR 391 \textit{Spatio-temporal Statistics for the Transition of Energy and Transport} (Project number 520388526) funded by the Deutsche Forschungsgemeinschaft (DFG, German Research Foundation).


\begin{thebibliography}{99}

\bibitem{Hoeffding1948}
Hoeffding, W. (1948).
A class of statistics with asymptotically normal distribution.
\emph{Annals of Mathematical Statistics} \textbf{19}, 293--325.

\bibitem{Neuhaus1977}
Neuhaus, G. (1977).
Functional limit theorems for $U$-statistics in the degenerate case.
\emph{Journal of Multivariate Analysis} \textbf{7}, 424--439.

\bibitem{ArconesGine1993}
Arcones, M. A. and Gin\'e, E. (1993).
Limit theorems for $U$-processes.
\emph{Annals of Probability} \textbf{21}, 1494--1542.

\bibitem{CCK2013}
Chernozhukov, V., Chetverikov, D. and Kato, K. (2013).
Gaussian approximations and multiplier bootstrap for maxima of sums of high-dimensional random vectors.
\emph{Annals of Statistics} \textbf{41}, 2786--2819.

\bibitem{CCK2015}
Chernozhukov, V., Chetverikov, D. and Kato, K. (2015).
Comparison and anti-concentration bounds for maxima of Gaussian random vectors.
\emph{Probability Theory and Related Fields} \textbf{162}, 47--70.

\bibitem{CCK2017}
Chernozhukov, V., Chetverikov, D. and Kato, K. (2017).
Central limit theorems and bootstrap in high dimensions.
\emph{Annals of Probability} \textbf{45}, 2309--2352.

\bibitem{CCKK2022}
Chernozhukov, V., Chetverikov, D., Kato, K. and Koike, Y. (2022).
Improved central limit theorem and bootstrap approximations in high dimensions.
\emph{Annals of Statistics} \textbf{50}, 2562--2586.

\bibitem{CCK2023}
Chernozhukov, V., Chetverikov, D. and Koike, Y. (2023).
Nearly optimal central limit theorem and bootstrap approximations in high dimensions.
\emph{Annals of Applied Probability} \textbf{33}, 2374--2425.

\bibitem{Chen2018}
Chen, X. (2018).
Gaussian and bootstrap approximations for high-dimensional $U$-statistics and their applications.
\emph{Annals of Statistics} \textbf{46}, 642--678.

\bibitem{ChenKato2019}
Chen, X. and Kato, K. (2019).
Randomized incomplete $U$-statistics in high dimensions.
\emph{Annals of Statistics} \textbf{47}, 3127--3156.

\bibitem{deJong1987}
de Jong, P. (1987).
A central limit theorem for generalized quadratic forms.
\emph{Probability Theory and Related Fields} \textbf{75}, 261--277.

\bibitem{DoblerPeccati2017}
D\"obler, C. and Peccati, G. (2017).
Quantitative de Jong theorems in any dimension.
\emph{Electronic Journal of Probability} \textbf{22}, Paper No. 2, 1--35.

\bibitem{DoblerPeccati2019}
D\"obler, C. and Peccati, G. (2019).
Quantitative CLTs for symmetric $U$-statistics using contractions.
\emph{Electronic Journal of Probability} \textbf{24}, Paper No. 5, 1--43.

\bibitem{NourdinPeccatiReinert2010}
Nourdin, I., Peccati, G. and Reinert, G. (2010).
Invariance principles for homogeneous sums: universality of Gaussian Wiener chaos.
\emph{Annals of Probability} \textbf{38}, 1947--1985.

\bibitem{Koike2019Wiener}
Koike, Y. (2019).
Gaussian approximation of maxima of Wiener functionals and its application to high-frequency data.
\emph{Annals of Statistics} \textbf{47}, 1663--1687.

\bibitem{Koike2021}
Koike, Y. (2021).
Notes on the dimension dependence in high-dimensional central limit theorems for hyperrectangles.
\emph{Japanese Journal of Statistics and Data Science} \textbf{4}, 257--297.
Corrected arXiv version: arXiv:1911.00160.

\bibitem{Koike2023}
Koike, Y. (2023).
High-dimensional central limit theorems for homogeneous sums.
\emph{Journal of Theoretical Probability} \textbf{36}, 1--45.

\bibitem{HuangEtAl2023}
Huang, K. H., Liu, X., Duncan, A. B. and Gandy, A. (2023).
A high-dimensional convergence theorem for $U$-statistics with applications to kernel-based testing.
\emph{Proceedings of Machine Learning Research} \textbf{195}, 3827--3918.

\bibitem[Imai and Koike(2025)]{ImaiKoike2025}
Imai, S. and Koike, Y. (2025).
Gaussian approximation for high-dimensional $U$-statistics with size-dependent kernels.
\emph{arXiv preprint} arXiv:2504.10866.

\bibitem{Vijaykumar2026}
Vijaykumar, S. (2026).
A scale-free density bound for Gaussian maxima.
\emph{arXiv preprint} arXiv:\allowbreak 2605.29066.

\bibitem{Pinelis2015}
Pinelis, I. (2015).
Rosenthal-type inequalities for martingales in 2-smooth Banach spaces.
\emph{Theory of Probability and Its Applications} \textbf{59}, 699--706.

\bibitem{Pinelis1994}
Pinelis, I. (1994).
Optimum bounds for the distributions of martingales in Banach spaces.
\emph{Annals of Probability} \textbf{22}, 1679--1706.

\bibitem{Talagrand2014}
Talagrand, M. (2014).
\emph{Upper and Lower Bounds for Stochastic Processes}.
Springer, Heidelberg.

\bibitem{Vershynin2018}
Vershynin, R. (2018).
\emph{High-Dimensional Probability: An Introduction with Applications in Data Science}.
Cambridge University Press, Cambridge.

\bibitem{KoltchinskiiLounici2017}
Koltchinskii, V. and Lounici, K. (2017).
Concentration inequalities and moment bounds for sample covariance operators.
\emph{Bernoulli} \textbf{23}, 110--133.

{
\bibitem[Carbery and Wright(2001)]{CarberyWright2001}
Carbery, A. and Wright, J. (2001).
Distributional and $L^q$ norm inequalities for polynomials over convex bodies
in $\R^n$.
\emph{Mathematical Research Letters} \textbf{8}, 233--248.
}

{
\bibitem[Decker et al.(2025)]{DeckerKongVolgushev2025}
Decker, C., Kong, D. and Volgushev, S. (2025).
Simultaneous hypothesis testing for comparing many functional means.
\emph{arXiv preprint} arXiv:2506.11889.
}

\bibitem[Arcones and Gin\'e(1992)]{ArconesGine1992}
Arcones, M. A. and Gin\'e, E. (1992).
On the bootstrap of $U$ and $V$ statistics.
\emph{The Annals of Statistics} \textbf{20}, 655--674.

\bibitem[Chwialkowski et al.(2014)]{Chwialkowski2014}
Chwialkowski, K., Sejdinovic, D. and Gretton, A. (2014).
A wild bootstrap for degenerate kernel tests.
\emph{Advances in Neural Information Processing Systems} \textbf{27}.

\bibitem[Leucht and Neumann(2013)]{LeuchtNeumann2013}
Leucht, A. and Neumann, M. H. (2013).
Dependent wild bootstrap for degenerate $U$- and $V$-statistics.
\emph{Journal of Multivariate Analysis} \textbf{117}, 257--280.


\bibitem[Cai and Jiang(2011)]{CaiJiang2011}
Cai, T. T. and Jiang, T. (2011).
Limiting laws of coherence of random matrices with applications to testing
covariance structure and construction of compressed sensing matrices.
\emph{Annals of Statistics} \textbf{39}, 1496--1525.

\bibitem{CaiLiuXia2013}
Cai, T. T., Liu, W. and Xia, Y. (2013).
Two-sample covariance matrix testing and support recovery in high-dimensional
and sparse settings.
\emph{Journal of the American Statistical Association} \textbf{108}, 265--277.

\bibitem[Karvonen and S\"arkk\"a(2019)]{KarvonenSarkka2019}
Karvonen, T. and S\"arkk\"a, S. (2019).
Gaussian kernel quadrature at scaled Gauss--Hermite nodes.
\textit{BIT Numerical Mathematics} \textbf{59}, 877--902.

\bibitem[Gretton et al.(2012)]{GrettonEtAl2012}
Gretton, A., Borgwardt, K. M., Rasch, M. J., Sch\"olkopf, B. and Smola, A. J.
(2012).
A kernel two-sample test.
\emph{Journal of Machine Learning Research} \textbf{13}, 723--773.

\bibitem[Gao and Shao(2023)]{GaoShao2023}
Gao, H. and Shao, X. (2023).
Two sample testing in high dimension via maximum mean discrepancy.
\emph{Journal of Machine Learning Research} \textbf{24}(304), 1--33.

\end{thebibliography}
\end{document}